\documentclass[12pt]{article}
\usepackage{amsfonts}

\usepackage{amsthm,amsmath,amssymb,amscd,verbatim,epsfig, graphicx}
\usepackage{graphicx}
\newtheorem{definition}{Definition}
\newtheorem{theorem}[definition]{Theorem}

\newtheorem{proposition}[definition]{Proposition}

\newtheorem{remark}[definition]{Remark}

\begin{document}

\title{Bifurcation of the ACT map }
\author{Bau-Sen Du\thanks{%
Institute of Mathematics, Academia Sinica, Taipei 115, TAIWAN,
dubs@math.sinica.edu.tw} , Ming-Chia Li\thanks{
Department of Applied Mathematics, National Chiao Tung University, Hsinchu
300, TAIWAN, mcli@math.nctu.edu.tw} \ and Mikhail Malkin\thanks{%
Department of Mathematics, Nizhny Novgorod State University, Nizhny
Novgorod, RUSSIA, malkin@uic.nnov.ru}}
\date{}
\maketitle

\begin{abstract}
In this paper, we study the Arneodo-Coullet-Tresser map $%
F(x,y,z)=(ax-b(y-z), bx+a(y-z), cx-dx^k+e z)$ where $a,b,c,d,e$ are real with $bd\neq 0$ and $k>1$ is an integer. We 
obtain
stability regions for fixed points of $F$ and symmetric period-2 points while $c$ and $e$ vary as parameters. 
Varying $a$ and $e$ as parameters, we show that there is a hyperbolic invariant set on which $F$ is conjugate to the 
full
shift on two or three symbols. We also show that chaotic behaviors of $F$ while $c$ and $d$ vary as parameters and 
$F$ is near an anti-integrable limit. Some numerical results indicates $F$ has Hopf bifurcation, strange
attractors, and nested structure of invariant tori.
\end{abstract}



\section{Introduction}

In this paper, we consider discrete dynamical systems induced by maps $F: {%
\mathbb{R}^3}\to {\mathbb{R}^3}$ of the form 
\begin{align}  \label{actmap}
F(x,y,z)=(ax-b(y-z), bx+a(y-z), cx-dx^k+e z),
\end{align}
where $a,b,c,d,e$ are real parameters with $bd\neq 0$ and $k>1$ is an
integer. These maps are introduced by Arneodo, Coullet and Tresser being motivated by the study of strange
attractors in a family of differential equations on ${\mathbb{R}^3}$ with homoclinic points of Shilnikov type,
refer to~\cite{Du1985}. Their numerical computations showed some interesting phenomena in dynamical behavior of
these maps (unpublished); in particular, for the one-parameter family of maps~\eqref{actmap} with $a=0.6, b=0.5, 
d=e=1,k=3
$ and $c$ as the parameter, they discovered (what may be regarded as) period doubling cascade and a strange
attractor. Rigorous results on characterization of periodic points of~\eqref{actmap}, as well as existence
of topological horseshoes for a region of parameters were obtained in~\cite{Du1985}. Nevertheless, since dynamical
behavior of maps~\eqref{actmap} is rich and multifarious in several regions of parameters, there remains many
interesting problems in rigorous understanding the dynamics and describing the bifurcation structure of~%
\eqref{actmap}; some of these problems will be considered in the present paper.

We will call the maps~\eqref{actmap} the Arneodo-Coullet-Tresser maps, or
ACT maps for short. Note that the Jacobian determinant of an ACT map $F$ is a constant; namely, $|\frac{%
\partial F(x,y,z)}{\partial (x,y,z)}|=(a^2+b^2)e$. Furthermore, if $e\neq 0$
then the map $F$ is a diffeomorphism with the inverse 
\begin{align*}
F^{-1}(x,y,z)=(\hat{x},\frac{-bx+ay}{a^2+b^2}+\hat{z}, \hat{z}),
\end{align*}
where $\displaystyle \hat{x}=\frac{ax+by}{a^2+b^2}$ and $\displaystyle \hat{z%
}=\frac{z-c\hat{x}+d\hat{x}^k}{e}$. So, $F$ is a polynomial automorphism of ${\mathbb{R}^3}$ and thus, in the
case of ACT maps, the Jacobian conjecture holds true (from this point of view, the ACT maps can be compared with
generalized H\'{e}non maps, which are also polynomial automorphisms, see~\cite{FriedlandMilnor1989, Qin2001} about 
the
history and results concerning the Jacobian conjecture, see~\cite{BassConnellWright1982, Rudin1995}).

When studying dynamical systems on noncompact manifolds, it is desirable to
be sure that the nonwandering set of the system is compact, in which case one might restrict oneself to this
nonwandering set in order to consider nontrivial dynamical behavior. In~\cite{LiMalkin2004}, a sufficient condition 
for the
nonwandering set of polynomial maps to be compact was obtained and an estimate for the size of the box containing 
the
nonwandering set was given. These results being applied to the case of ACT maps imply (see~\cite[Proposition 15%
]{LiMalkin2004}) that the nonwandering set $\Omega(F)$ (as well as the set of bounded orbits) of the ACT map $F$
lies in the box 
\begin{align*}
\left\{(x,y,z)\in {\mathbb{R}^3}: |x|\leq M, |y|\leq \frac{a^2+b^2+2|a|+1}{%
|b|}M , |z|\leq \frac{a^2+b^2+|a|}{|b|}M\right\},
\end{align*}
where 
\begin{align*}
M=\sqrt[k-1]{\frac{|a^2e+b^2e|+|a^2+b^2-bc+2ae|+|2a+e|+1}{|bd|}}.
\end{align*}
Note that the above result of the nonwandering set is a generalization of a
result in~\cite{Du1985} which asserts that the set of periodic points of the ACT map is contained in the box.

We will discuss (see Sections 3 and 4) chaotic dynamical behavior of the ACT
maps in some parameter region. Usually chaotic behavior of a dynamical system is associated with positive
topological entropy. For maps not uniformly continuous with non-compact domain one needs to be very specific to 
define
the topological entropy $h_{top}(F)$ to be $h_{top}(F|_{\Omega(F)})$. On the other hand, since $F$ is a polynomial
automorphism of ${\mathbb{R}^3}$, it can be easily extended to the homeomorphism $\bar{F}$ of the canonical one 
point
compactification $\bar{\mathbb{R}}^3:={\mathbb{R}^3}\bigcup\{\infty\}$ with $%
\bar{F}(\infty)=\infty$ and so one can also define $h_{top}(F)$ to be $h_{top}(\bar{F})$ and the two definitions of 
the topological entropy agree
because $h_{top}(\bar{F})=h_{top}(\bar{F}|_{\Omega(\bar{F})})$ and $\infty$ is an
isolated nonwandering point for $\bar{F}$.

The paper is organized as follows. In Section 2, we obtain stability regions
for fixed points of the ACT maps and symmetric period-2 points for two-parameter family of the ACT maps with $c$
and $e$ as the parameters. These stability regions are used later in Section 5 to explain period doubling
bifurcation and discrete Andronov-Hopf bifurcations. We also observe that the intersection of these stability
regions is a point at which the Jacobiant determinant of the ACT map takes value $1$ (and so the system becomes
conservative) and which gives rise to a nested family of invariant two-dimensional tori, and as numerical study 
shows, this
nested structure of invariant tori persist for some interval of parameter $c$ with parameter $e$ fixed (so that
that system remains conservative). In Section 3, we study chaotic behavior of the ACT maps for two-parameter
family with $a$ and $e$ as parameters. We prove that for $a$ and $e$ small there is a hyperbolic invariant set which
is conjugate to the full shift on two or three symbols depending on the evenness of the number $k$ 
in~\eqref{actmap}.
In Section 4, based on a result in~\cite{LiMalkin2004diff}, we study the
chaotic dynamics of the ACT maps for two-parameter family with $c$ and $d$ as parameters. While $c\to \infty$ and 
$d\to\infty$ in such a way that $\frac{d}{c}=\text{constant}>0$. Let $\lambda=\frac{1}{c}$. Then for all $%
|\lambda|$ sufficiently small, the ACT map $F=F_\lambda$ has a closed invariant set $\Lambda_\lambda$ such that $%
F_\lambda|\Lambda_\lambda$ is conjugate to the full shift on either two or three symbols depending on whether $k$ is
even or odd, respectively. In Section 5, we discuss numerical results on
bifurcations of the ACT maps including Hopf bifurcations and resonances, appearance and changing of invariant 
circles, strange
attractors, nested structure of invariant tori. In the Appendix, we derive the stability criterion for $3\times 3$ 
real
matrices, which is used in Section 2 (and is more convenient for our purposes than the Shur-Cohn criterion).

While studying parameter-dependence property of $F$, we will denote $F$ as $%
F_c$ (resp. $F_e$) to stress that $c$ varies (resp. $e$ varies) as other parameters are fixed. Similarly, we
denote $F_{c,e}$ when $c$ and $e$ are the only parameters that vary.

\section{Stability regions for fixed points and symmetric period-two points}

It is clear that the origin is a fixed point of $F$. By solving $%
F(x,y,z)=(x,y,z)$, we obtain that for even $k$, the map $F$ has a unique nontrivial fixed point at 
\begin{align*}
p_1=(x_1,\frac{a^2+b^2-a}{b}x_1,\frac{(a-1)^2+b^2}{b}x_1)
\end{align*}
where $x_1=\sqrt[k-1]{\frac{bc-(1-e)[(a-1)^2+b^2]}{bd}}$, and that for the
case when $k$ is odd and $\frac{bc-(1-e)[(a-1)^2+b^2]}{bd}>0$, the map $F$ has exactly two nontrivial
fixed points at $\pm p_1$. Moreover, by solving the system of equations $F(x,y,z)=(-x,-y,-z)$ and $%
F(-x,-y,-z)=(x,y,z)$, we get that $F$ has periodic points of period $2$ symmetric to the origin, say $\pm p_2$, if 
and only if $%
k$ is odd and $\frac{bc-(1+e)[(a+1)^2+b^2]}{bd}>0$; under these conditions $p_2$ is given
by 
\begin{align*}
p_{2}=(x_2,\frac{-a^2-b^2-a}{b}x_2,\frac{-(a+1)^2-b^2}{b}x_2),
\end{align*}
with $x_2=\sqrt[k-1]{\frac{bc-(1+e)[(a+1)^2+b^2]}{bd}}$. We call $\pm p_2$
(when exist) the \emph{symmetric period-two points} of $F$.

In next section, we will be concerned with stability regions for fixed
points and period-$2$ points of the ACT maps, i.e., regions in the parameter space for which these points are 
stable. More
precisely, we use the following definitions.

\begin{definition}
Let $\mathbf{x}\mapsto G_{\mathbf{v}}(\mathbf{x}), \mathbf{x}\in X\subset {%
\mathbb{R}^m}, \mathbf{v}\in V\subset{\mathbb{R}^l}$ be a family of $C^1$ maps with parameter $\mathbf{v}$, and let 
$\mathbf{x}^*_{%
\mathbf{v}}$ be a fixed point of $G_{\mathbf{v}}$ for each $\mathbf{v}\in V$. A subset $J\subset V$ is called the 
\emph{stability region%
} for the family of fixed points $\mathbf{x}^*_{\mathbf{v}}$ if for each $\mathbf{v}\in J$, $G_{\mathbf{v}}$ is 
\emph{%
linearly stable} at $\mathbf{x}^*_{\mathbf{v}}$ and for each $\mathbf{v}\notin J$, $G_{\mathbf{v}}$ is not 
\emph{linearly stable} at $%
\mathbf{x}^*_{\mathbf{v}}$, i.e., for each $\mathbf{v}\in J$ all eigenvalues of the Jacobian matrix 
$\frac{\partial}{\partial \mathbf{x}}%
G_{\mathbf{v}}(\mathbf{x}^*_{\mathbf{v}})$ lie in $\{z\in {\mathbb{C}}: |z|<1\}$ and for each $\mathbf{v}\notin J$ 
there is an
eigenvalue of the Jacobian matrix $\frac{\partial}{\partial \mathbf{x}}G_{\mathbf{v}}(\mathbf{x}^*_{\mathbf{v}%
})$ in $\{z\in {\mathbb{C}}: |z|\geq 1\}$.
\end{definition}

Similarly, we can define a \emph{stability region} for a family of a
periodic orbit of period $n$ by replacing $G_{\mathbf{v}}$ by $G^n_{\mathbf{v}}$.

For the ACT family $F$, the Jacobian matrix of $F$ at a point $(x,y,z)$ is 
\begin{align*}
\frac{\partial F(x,y,z)}{\partial (x,y,z)}= 
\begin{bmatrix}
a & -b & b \\ 
b & a & -a \\ 
c-kdx^{k-1} & 0 & e%
\end{bmatrix}%
\end{align*}
and its characteristic polynomial is 
\begin{equation*}
P(\lambda)=\lambda^3-(2a+e)\lambda^2+[a^2+b^2+2ae-bc+kbdx^{k-1})]%
\lambda-(a^2+b^2)e.
\end{equation*}

Note that the determinant of the Jacobian matrix of $F$ is constant, that
is, $\displaystyle |\frac{\partial F(x,y,z)}{\partial (x,y,z)}|=(a^2+b^2)e$,
and if $e\neq 0$ then the map $F:{\mathbb{R}^3}\to {\mathbb{R}^3}$ is a diffeomorphism with the inverse 
\begin{equation*}
F^{-1}(x,y,z)=(\hat{x}, \frac{-b x+a y}{a^2+b^2}+\hat{z}, \hat{z}), 
\end{equation*}
where $\hat{x}=\displaystyle \frac{a x+b y}{a^2+b^2}$ and $\displaystyle 
\hat{z}=\frac{z-c\hat{x}+d\hat{x}^k}{e}$. So, $F$ is a polynomial
automorphism of ${\mathbb{R}^3}$ and thus, in the case of ACT maps, the Jacobian conjecture
holds true (about the history and results concerning the Jacobian conjecture, see \cite{BassConnellWright1982}, 
\cite%
{Rudin1995}).

To determine the local stability of the fixed points, we need to know
whether eigenvalues of the Jacobian matrix of ACT map lies inside the unit circle in the complex plane, i.e., 
whether
roots of the characteristic polynomial of this matrix lies inside the unit circle. We will use the following
definition.

\begin{definition}
We say that a polynomial $P$ over ${\mathbb{R}}$ is \emph{stable} if all of
its roots lie in $\{z\in {\mathbb{C}}: |z|<1\}$.
\end{definition}

Note that polynomials of this type of stability are called sometimes \emph{%
discrete stable} or \emph{Shur stable} to distinguish from the situation, when all roots of polynomial have negative
real parts (this applies usually to characteristic polynomials of linearization of vector fields at stable fixed
points), in which case the term \emph{Hurwitz stability} is sometimes used. There are several criteria to
determine stability of a polynomial in terms of its coefficients (e.g., the
Shur-Cohn criterion, see \cite{Marden1989}). For our purposes (taking in mind
bifurcation parameters of ACT maps) it will be convenient to use the following one, which is derived in the 
Appendix.

\begin{proposition}
\label{polystable} A polynomial $P(\lambda)=\lambda^3+A\lambda^2+B\lambda+D%
\in \mathbb{R}[\lambda]$ is stable if and only if 
\begin{align*}
\text{$|D|<1$ and $\max\{-P(1),P(-1)\}<0<\hat{\alpha}$,}
\end{align*}
where $\hat{\alpha}=\hat{\alpha}(P):=-D^2+AD-B+1$.
\end{proposition}

We apply now the above criterion to a one-parameter family of polynomials,
which are the characteristic polynomials of the Jacobian matrix of ACT maps.
The proof of the following proposition is postponed to the Appendix.

\begin{proposition}
\label{bifstable} Let $P_\beta(\lambda)=\lambda^3+A\lambda^2+(B+\beta)%
\lambda+D$ be a one-parameter family of polynomials in $\mathbb{R}[\lambda]$ with parameter $\beta\in {\mathbb{R}}$. 
Then

\begin{enumerate}
\item The following three statements are equivalent:

\begin{enumerate}
\item $|D|<1$ and $|A-D|<2$.

\item There exists a $\beta$ such that all roots of $P_\beta$ lie in $\{z\in%
\mathbb{C}:|z|<1\}$.

\item There is a unique open interval $I\subset \mathbb{R}$, which will be
called the \emph{stable interval} of the family $P_\beta$, such that for each $\beta\in I$, $P_\beta$ is stable and
for each $\beta\notin I$ there is a root of $P_\beta$ in $\{z\in\mathbb{C}:|z|\geq 1\}.$
\end{enumerate}

\item The stable interval of $P_\beta$ $($when exists, see item $1)$ is
equal to $(\alpha, \hat{\alpha}_0)$, where $\alpha=\max\{-P_0(1),P_0(-1)\}$ and $\hat{\alpha}_0=-D^2+AD-B+1$. 
Moreover,
for $\beta=\hat{\alpha}_0$ the polynomial $P_\beta$ has two complex conjugate roots in $\{z\in \mathbb{C}:
|z|=1\}\backslash \{-1,1\}$.
\end{enumerate}
\end{proposition}

Now we are in position to study the stability regions for the map $F$. In
the following Theorems~\ref{stregorigin}-\ref{stregsym}, we deal with stability regions $%
J_{\text{tr}}(F_c)$ and $J_{\text{tr}}(F_{c,e})$ (recall that the subscripts here indicate the only
parameters that vary). Note that if we find some functions $f_1(e)$, $f_2(e)$ of variables $e$, such that 
$J_{\text{%
tr}}(F_{c,e})=\{(e,c)\in {\mathbb{R}^2}: f_1(e)<c<f_2(e)\}$ then we will have for any $c$ with $J_{\text{tr}%
}(F_c)\neq \o $, that $J_{\text{tr}}(F_c)=\{c\in {\mathbb{R}}: f_1(e)<c<f_2(e)\}$, i.e., $J_{\text{tr}}(F_c)$ is 
described by
the same inequalities. So in this case it is enough to give the corresponding formulas for $J_{\text{tr}}(F_{c,e})$
only.

\begin{theorem}[stability regions for the trivial fixed point]
\label{stregorigin} \ Let $F$ be the ACT family with $b\neq 0$. Let $J_{%
\text{tr}}(F_c)$ $($resp. $J_{\text{tr}}(F_{c,e}))$ denote the stability region of the origin for $F_c$ $($resp. for 
$F_{c,e}$$)$. Then

\begin{enumerate}
\item $J_{\text{tr}}(F_c)\neq \o $ if and only if 
\begin{align}
\text{ $-1<(a^2+b^2)e<1$ \hspace{1ex} and \hspace{1ex} $%
2a-2<(a^2+b^2-1)e<2a+2.$}
\end{align}

\item For $F_{c,e}$, the following two statements hold:

\begin{enumerate}
\item Suppose $a^2+b^2-1\leq 0$, then $J_{\text{tr}}(F_{c,e})\neq \o $.

\item Suppose $a^2+b^2-1>0$ , then $J_{\text{tr}}(F_{c,e})\neq \o $ if and
only if 
\begin{align}
\displaystyle \max\left\{\frac{2a-2}{a^2+b^2-1}, \frac{-2a-2}{a^2+b^2-1}%
\right\}<\frac{1}{a^2+b^2}.
\end{align}
\end{enumerate}

\item If $J_{\text{tr}}(F_{c,e})\neq \o $, then 
\begin{equation*}
J_{\text{tr}}(F_{c,e})=\{(e,c)\in{\mathbb{R}^2}: \max\{-c_1(e),c_{-1}(e)\}<-bc<\hat{c}(e)\},
\end{equation*}
where 
\begin{align}
\begin{array}{l}
c_1(e)={(1-e)[(a-1)^2+b^2]},\quad c_{-1}(e)=\displaystyle{-(1+e)[(a+1)^2+b^2]%
} \\ 
\ \quad \mbox{and}\quad \hat{c}(e)=\displaystyle{%
-(a^2+b^2-1)[(ae-1)^2+b^2e^2]}.%
\end{array}%
\end{align}
$($See Figure 1 and the remarks below it. $)$
\end{enumerate}
\end{theorem}

\begin{proof} Since $b\neq 0$, we may use $\beta=-bc$ as parameter. 
Then the characteristic polynomial of the Jacobian matrix of $F$ at the origin is  
$$P_\beta(\lambda)=\lambda^3-(2a+e)\lambda^2+(a^2+b^2+2ae+\beta)\lambda-(a^2+b^2)e.$$ 
By applying item $1$ of Proposition~\ref{bifstable} to the above family $P_\beta$ with parameter $\beta$, we obtain 
that equation (4) is equivalent to that the stable interval of $P_\beta$ is not empty  and so is equivalent to that 
the stability region of $F$ is not empty.  Item 1 is completed.  

The inequalities $a^2+b^2-1\leq 0$ and $b\neq 0$ imply $|a|<1$, so the numbers $2a-2$ and $2a+2$ are of opposite 
signs, and therefore one gets item $2 (a)$ from item $1$ by taking $e=0$.

For item 2(b), let us denote $e_\pm=\displaystyle \frac{\pm 1}{a^2+b^2}$, $e_l=\displaystyle \frac{2a-2}{a^2+b^2-1}$ 
and $e_r=\displaystyle \frac{2a+2}{a^2+b^2-1}$.
The existence of $e$ satisfying condition (4) is equivalent to the fact that the two intervals $(e_-,e_+)$ and 
$(e_l,e_r)$ overlap, i.e., $e_l<e_+$ and $e_-<e_r$.  It is easy to see that the last two inequalities are the same 
as (5). 

The stable interval of $P_\beta$ is given by $\max\{-P_0(1),P_0(-1)\}<\beta<\hat{\alpha}_0$, (for the definition of 
$\hat{\alpha}_0$ see item $2$ of Proposition~\ref{bifstable}).  By evaluating the values $P_0(1)$, $P_0(-1)$ and 
$\hat{\alpha}_0$ and using the fact that $\beta=-bc$, item $3$ is proved. 
\end{proof}

\begin{figure}[h]
\caption{The graphs of the bifurcation curves $c=\frac{c_1(e)}{b}$, $c=\frac{%
c_{-1}(e)}{-b}$ and $c=\frac{\hat{c}(e)}{-b}$, indicated simply as $c_1,
c_{-1}$ and $\hat{c}$, are shown in the $(e,c)$-plane $("e"$ is the
horizontal axis and $"c"$ is the vertical one$)$. The dashed lines in the
figures are $e=\frac{\pm 1}{a^2+b^2}$. In figures $1(i)$-$(iii)$, the
stability regions $J_{\text{tr}}(F_{c,e})$ (shaded in black) are shown for
three cases: when $a^2+b^2-1=0,>0$ and $<0$, all together with $b<0$, namely 
$(i)$ $a=0.6$ and $b=-0.8$, $(ii)$ $a=0.2$ and $b=-1.4$, and $(iii)$ $a=0.1$
and $b=-0.8$. Figure $1(iv)$ corresponds to the subcase of $(ii)$ when the
stability region has two sides $(M^{\prime\prime}$ lies inside the strip
between the dashed lines$)$; here $a=0.85$ and $b=-1$. }%
\centerline{\hspace{0ex}\epsfig{file=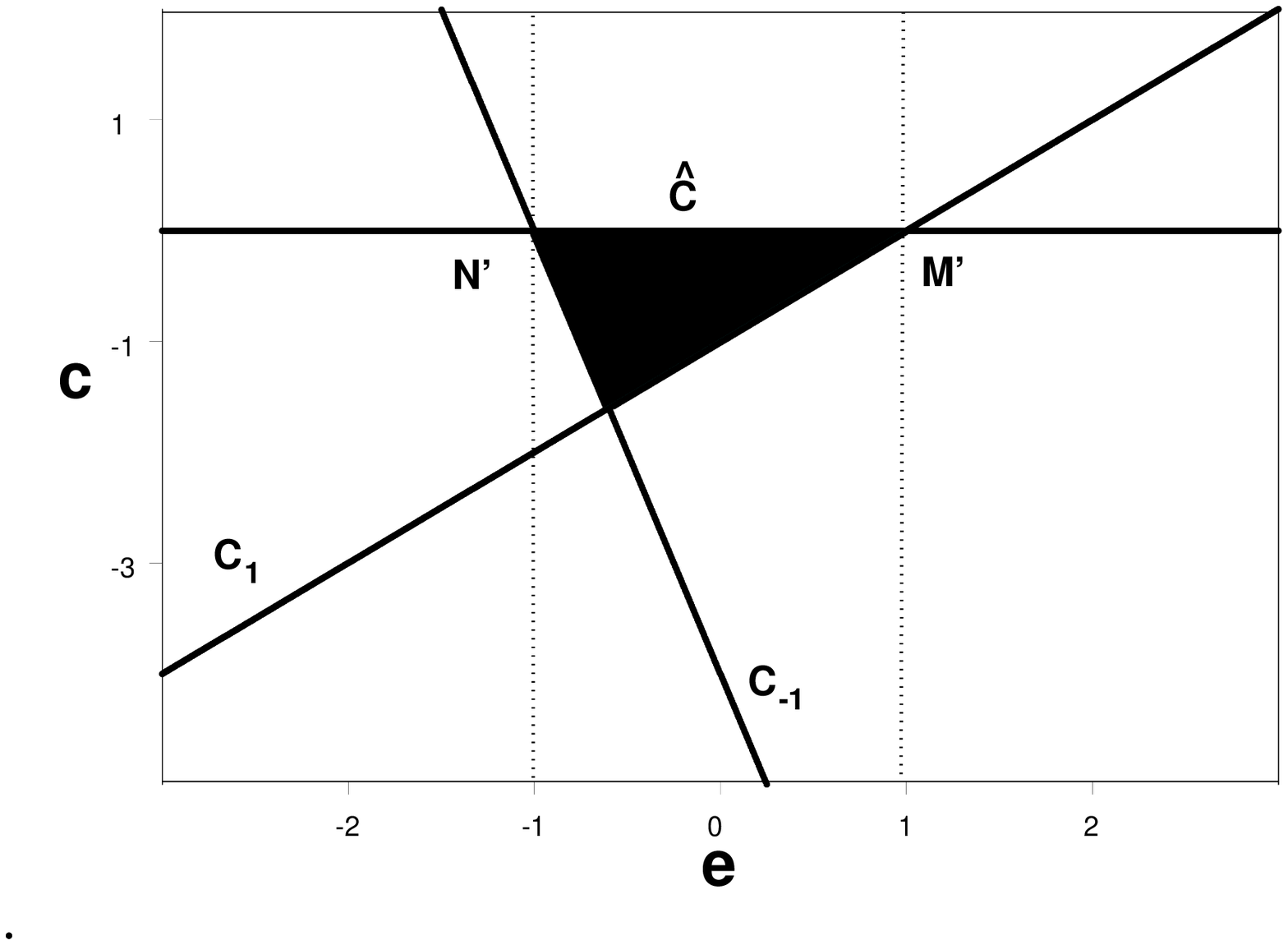,width=5.7cm,height=4.5cm}
\hspace{4ex}\epsfig{file=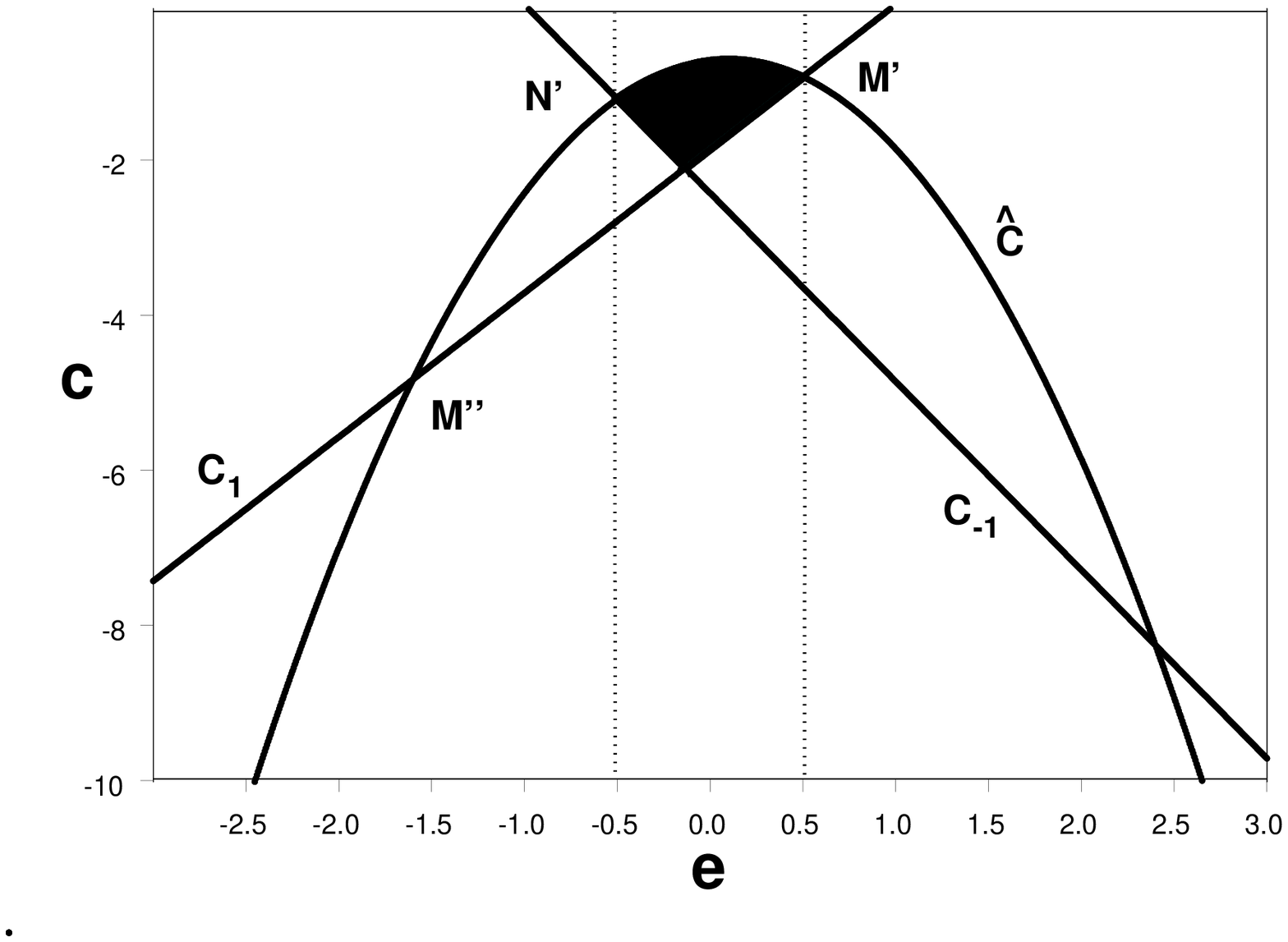,width=5.7cm,height=4.5cm}
} \centerline{\ \ \ $(i)$\hspace{6.1cm}$(ii)$} 
\centerline{\hspace{0ex}\epsfig{file=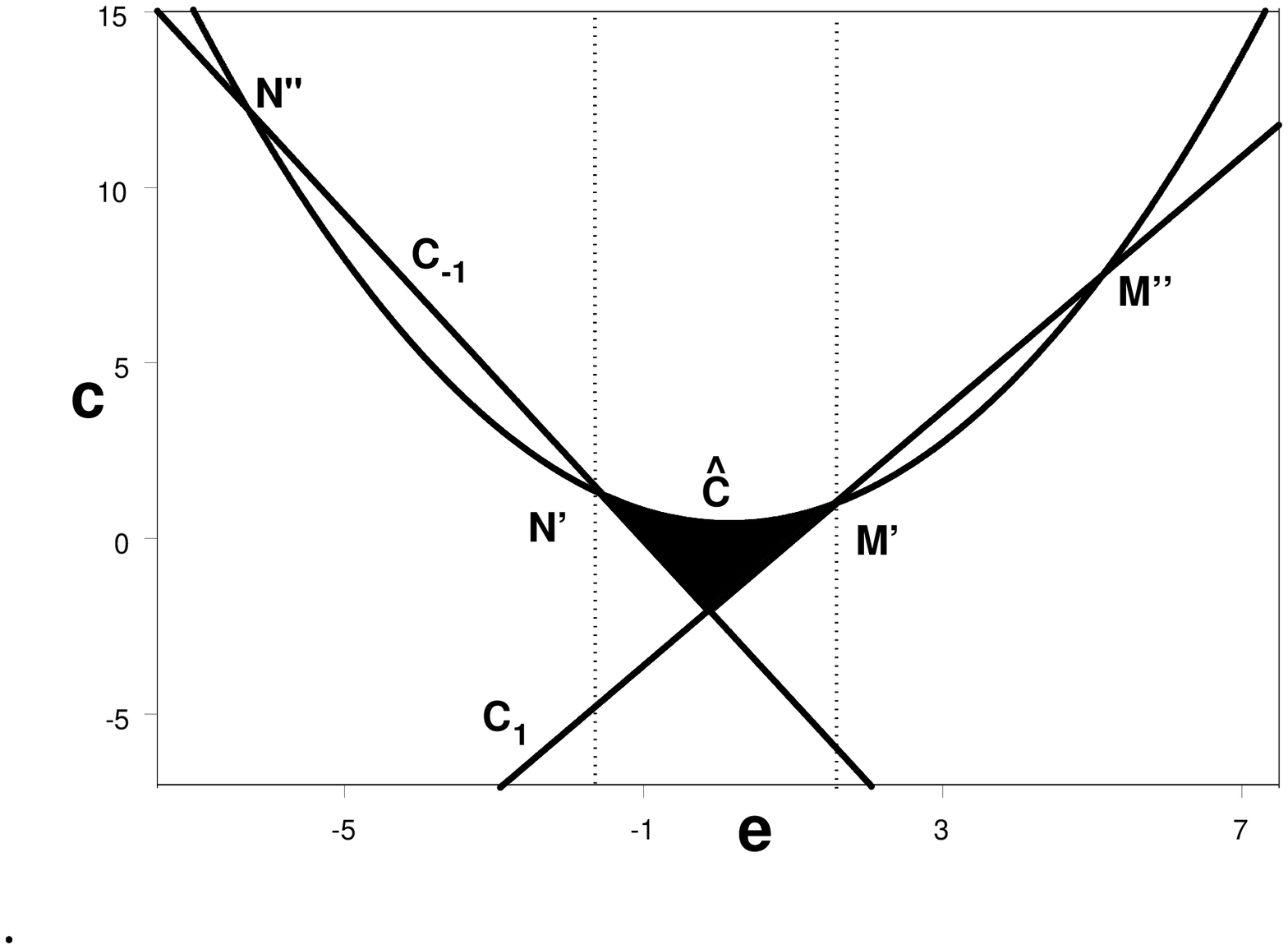,width=5.7cm,height=4.5cm}
\hspace{4ex}\epsfig{file=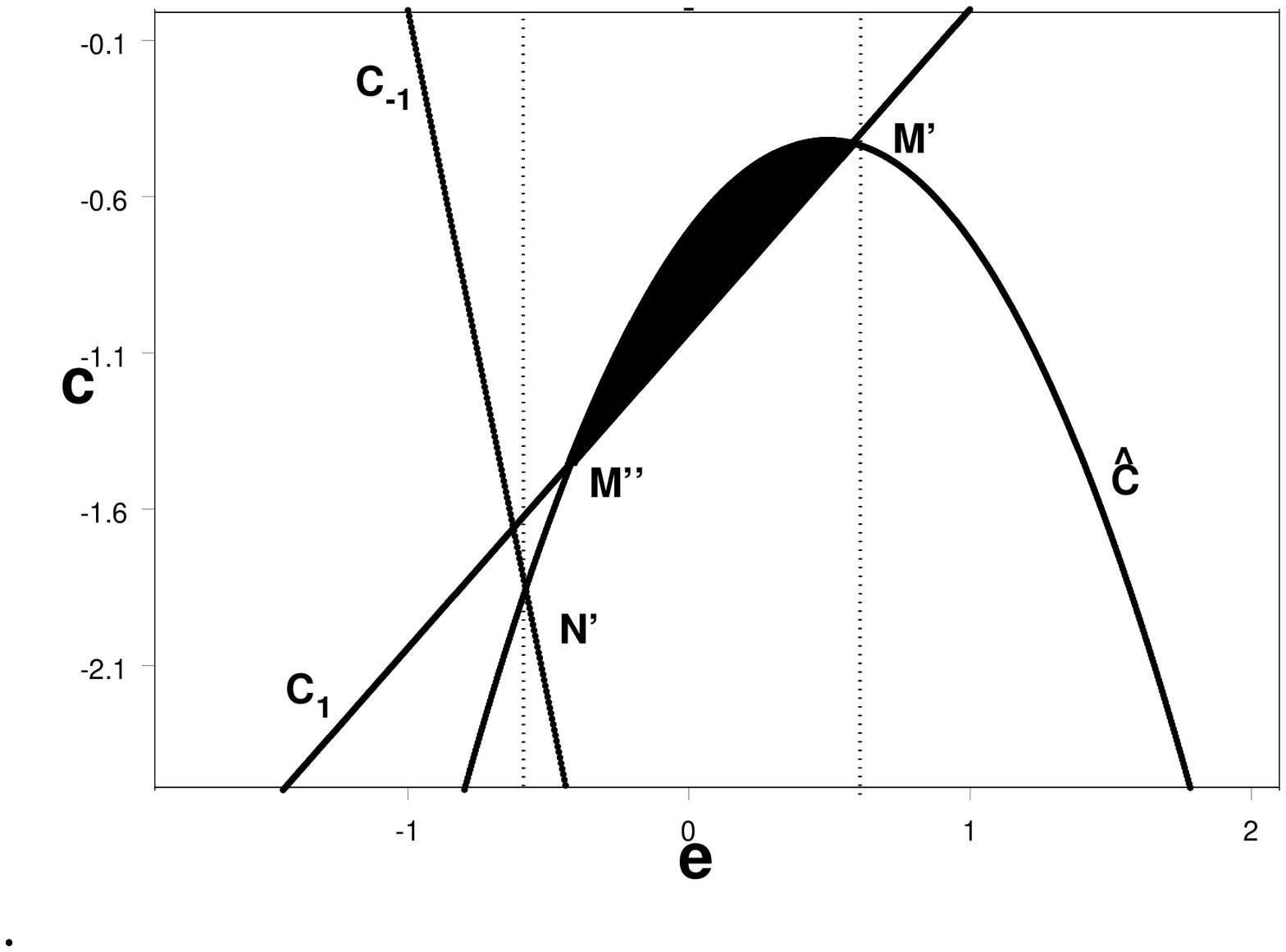,width=5.7cm,height=4.5cm}} \centerline{\ \
\ $(iii)$\hspace{5.9cm}$(iv)$}
\par
\qquad Let us give some remarks on Figure 1. The dashed lines there are $%
e=\pm \frac{1}{a^2+b^2}$, which corresponds to the cases when $D=\pm 1$ in Propositions $1$ and $2$. So the
stability region $J_{\text{tr}}(F_{c,e})$ must belong to the strip between the dashed lines.
It is easy to see that the intersection of the bifurcation curves $c=\frac{\hat{c}(e)}{-b}$ and $c=\frac{c_1(e)}{b}$
consists of either one or two points depending on whether $a^2+b^2-1$ is zero or not. In the former case, the
point of intersection has coordinates $(e,c)=(\frac{1}{a^2+b^2},0)$, which corresponds to the eigenvalues $%
\lambda_1=1$ and $\lambda_{2,3}=a\pm i\sqrt{1-a^2}$. In the latter case, the two intersection points are $%
M^\prime(e^\prime, c^\prime)$ and $M^{\prime\prime}(e^{\prime\prime}, c^{\prime\prime})$, where 
$c^\prime=\frac{(1-\frac{1}{a^2+b^2})[(a-1)^2+b^2]}{b}$, $e^\prime=\frac{1}{%
a^2+b^2}$, $c^{\prime\prime}=\frac{(1-\frac{2a-2}{a^2+b^2-1})[(a-1)^2+b^2]}{b}$, and $%
e^{\prime\prime}=\frac{2a-2}{a^2+b^2-1}$, which corresponds to the eigenvalues $\lambda_1^\prime=1$, $%
|\lambda_2^\prime|=|\lambda_3^\prime|=1$ and $\lambda_1^{\prime\prime}=\lambda_2^{\prime\prime}=1$, 
$\lambda_3^{\prime%
\prime}\in{\mathbb{R}}$. Note that the coordinates of $M^\prime$ and $M^{\prime\prime}$ satisfy the equalities $D=1$
and $D=A+2$ respectively (see conditions in Proposition~\ref{polystable}, item 1(a)). Furthermore, the
condition (5) on existence of stability region $J_{\text{tr}}(F_{c,e})$ implies that $e^{\prime\prime}<e^\prime$ if 
$%
a^2+b^2-1>0$. Similar geometric interpretation of the conditions of
Proposition~\ref{bifstable} and Theorem~\ref{stregorigin} can be done in terms of intersections, say $N^\prime$ and 
$N^{\prime\prime}$, of
the lines $c=\hat{c}$ and $c=c_{-1}$. Note that $J_{\text{tr}}(F_{c,e})$ has either two "sides" or three "sides"
depending on whether both points $M^{\prime\prime}$ and $N^{\prime\prime}$ lie outside the strip $|e|<\frac{1%
}{a^2+b^2}$ or one of $M^{\prime\prime}$ and $N^{\prime\prime}$ lies inside $($one can easily see that the points $%
M^{\prime\prime}$ and $N^{\prime\prime}$ cannot lie simultaneously inside the strip$)$.
\end{figure}

Next we study the stability region of the nontrivial fixed point(s), say $p_1
$ for even $k$ and $\pm p_1$ for odd $k$. Let $x_1$ be the $x$-coordinate of $p_1$, then one easily gets that $%
x_1=\sqrt[k-1]{\frac{bc-c_1}{bd}}$, where $c_1=(1-e)[(a-1)^2+b^2]$ (the same as in item 3 of Theorem~\ref%
{stregorigin}). Similar to Theorem~\ref{stregorigin}, applying Proposition~\ref{bifstable} to the
characteristic polynomial of the Jacobian matrix of $F$ at the nontrivial fixed point(s) $p_1$ (or $\pm p_1$), which
is 
\begin{equation*}
P_\beta(\lambda)=\lambda^3-(2a+e)\lambda^2+(a^2+b^2+2ae+\beta)%
\lambda-(a^2+b^2)e, 
\end{equation*}
where $\beta=-bc+kbdx_1^{k-1}$, we have the following result.

\begin{theorem}[stability regions for the nontrivial fixed point(s)]
\ Let $J_{\text{non}}(F_c)$ $($resp. $J_{\text{non}}(F_{c,e})$$)$ denote the
stability region of the nontrivial fixed point$($s$)$ for the map $F_c$ $($resp. for $(F_{c,e})$$)$. Then

\begin{enumerate}
\item If $k$ is odd and $bd<0$ then $J_{\text{non}}(F_c)=J_{\text{non}%
}(F_{c,e})=\o $. If $k$ is even or $bd>0$, then $J_{\text{non}}(F_c)\neq \o $
if and only if $J_{\text{tr}}(F_c)\neq \o $; moreover $J_{\text{non}}(F_{c,e})\neq \o $ if and only if 
$J_{\text{tr}}(F_{c,e})\neq 
\o $.

\item If $J_{\text{non}}(F_{c,e})\neq \o $, then  
\begin{align*}
J_{\text{non}}(F_{c,e})=&\Big\{(e,c)\in{\mathbb{R}^2}: \\
&\quad -\frac{kc_1(e)+\hat{c}(e)}{k-1}<-bc<\min\{-c_1(e),-\frac{%
kc_1(e)+c_{-1}(e)}{k-1}\}\Big\},
\end{align*}
where $\hat{c}(e)$, $c_{1}(e)$ and $c_{-1}(e)$ are in item 3 of Theorem~\ref%
{stregorigin}.

\item For $F_{c,e}$, the point $(-c_1(e),e)$ is on the boundary of the
stability region of the nontrivial fixed point(s) if and only if it is on the boundary of the stability region of 
the
origin.

See Figure~2 for the regions of $J_{\text{non}}(F_{c,e})$ aside with $J_{%
\text{tr}}(F_{c,e})$.
\end{enumerate}
\end{theorem}

\begin{figure}[h]
\centerline{\hspace{0ex}\epsfig{file=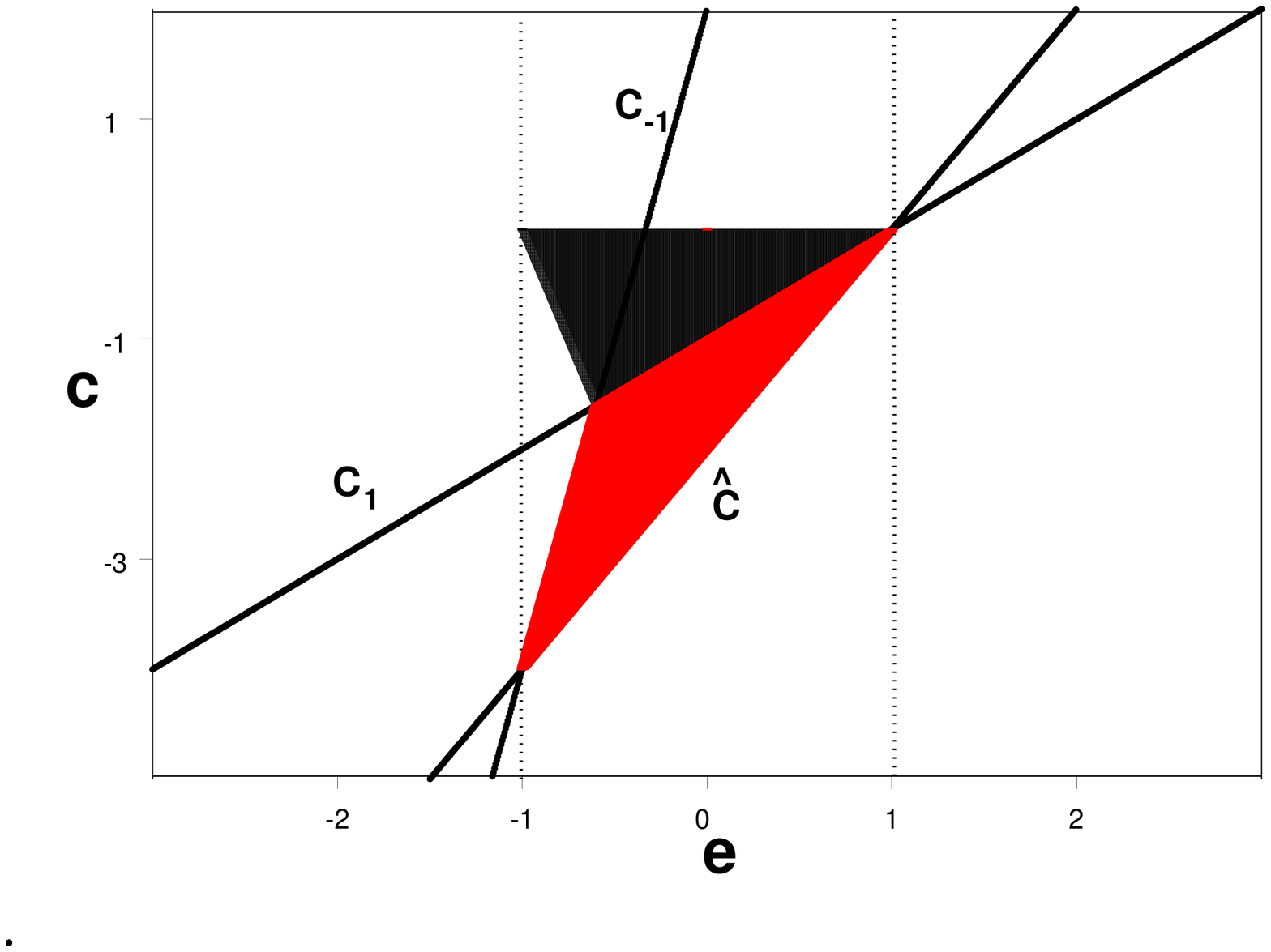,width=5.7cm,height=5.7cm}
\hspace{4ex}\epsfig{file=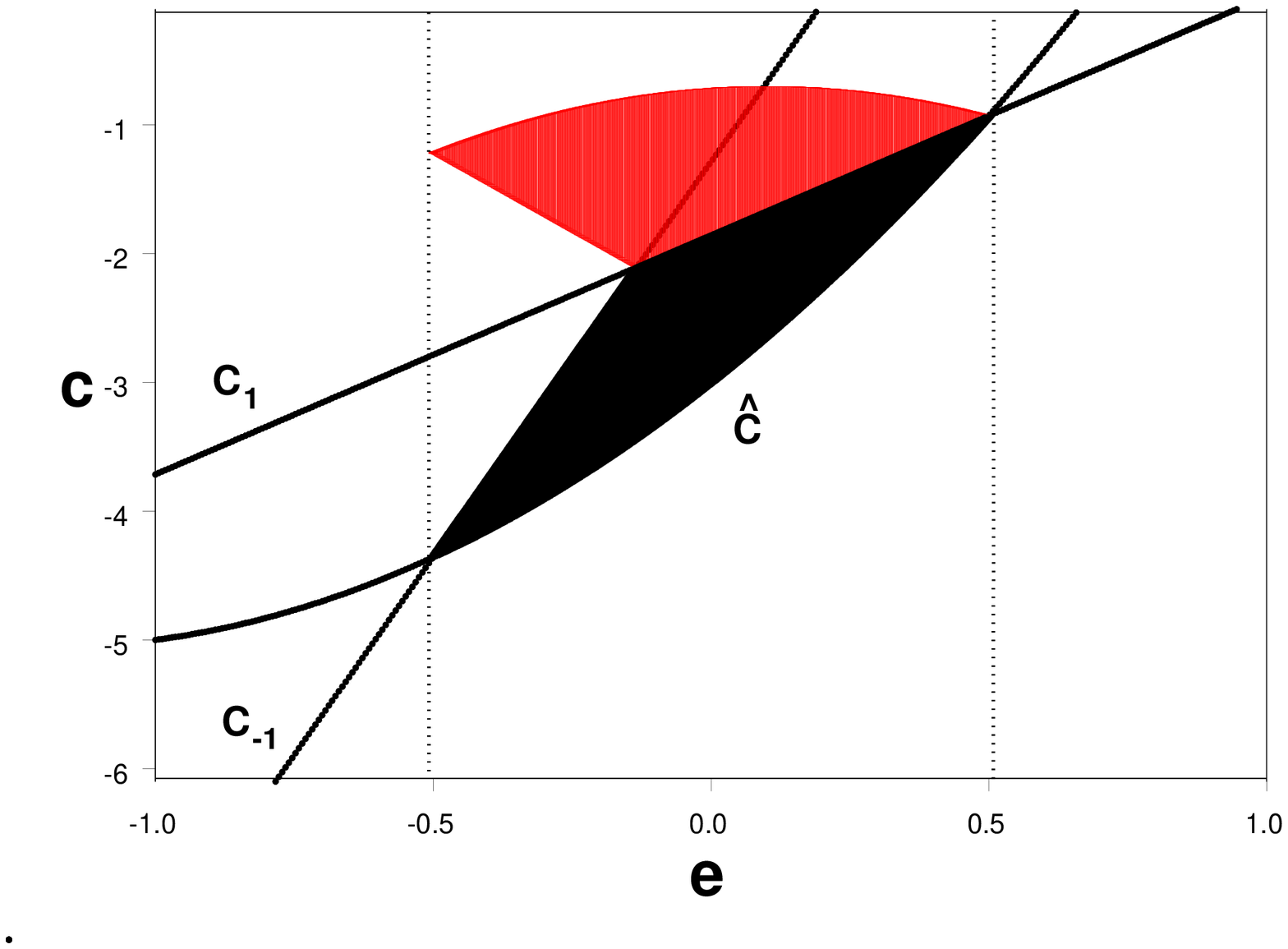,width=5.7cm,height=5.7cm}} \centerline{\ \
\ $(i)$\hspace{6.1cm}$(ii)$} 
\centerline{\hspace{0ex}\epsfig{file=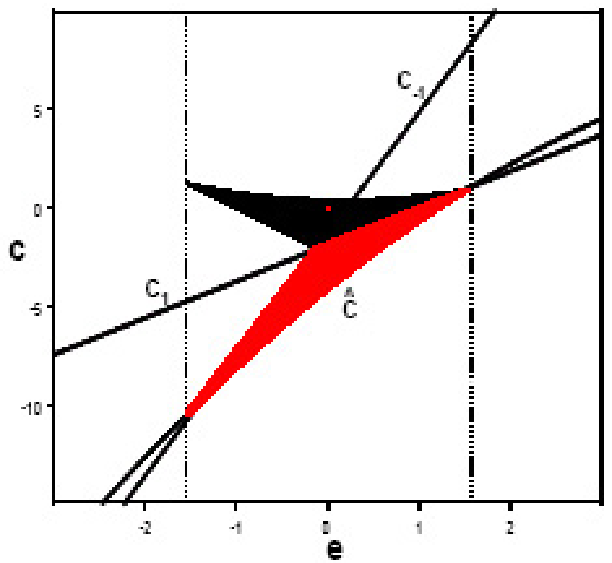,width=5.7cm,height=5.7cm}
} \centerline{\ \ \ $(iii)$}
\caption{The graphs of the bifurcation curves $c=\frac{kc_1(e)+\hat{c}(e)}{%
b(k-1)}$, $c=\frac{c_{1}(e)}{b}$ and $c=\frac{kc_1(e)+c_{-1}(e)}{b(k-1)}$,
indicated simply as $\hat{c}, c_{1}$ and $c_{-1}$, are shown in the $(e,c)$%
-plane $("e"$ is the horizontal axis and $"c"$ is the vertical one$)$. The
dashed lines in the figures are $e=\frac{\pm 1}{a^2+b^2}$. In figures $1(i)$-%
$(iii)$, the stability regions $J_{\text{non}}(F_{c,e})$ (shaded in black)
and $J_{\text{tr}}(F_{c,e})$ (shaded in gray) are shown for three cases:
when $a^2+b^2-1=0,>0$ and $<0$, all together with $b<0$, namely $(i)$ $a=0.6$
and $b=-0.8$, $(ii)$ $a=0.2$ and $b=-1.4$, and $(iii)$ $a=0.1$ and $b=-0.8$. 
}
\end{figure}

For the symmetric period-two points at $\pm p_2$ of the ACT family $F$ with
odd $k$, we have 
\begin{equation*}
DF^2(p_2)=DF(p_2)\cdot DF(-p_2)=[DF(p_2)]^2
\end{equation*}
and so the eigenvalues of $DF^2(p_2)$ are the square of the eigenvalues of $DF(p_2)$. It follows that the stability 
of $\pm p_2$ is
determined by whether the characteristic polynomial of $DF(p_2)$ has all roots in the open unit disk. Let $x_2$ be
the $x$-coordinate of $p_2$, then $x_2=\sqrt[k-1]{\frac{bc+c_{-1}}{bd}}$, where $c_{-1}$ is in item $3$ of
Theorem~\ref{stregorigin}. Applying Proposition~\ref{bifstable} to the characteristic polynomial of $DF(p_2)$,
which is 
\begin{equation*}
P_\beta(\lambda)=\lambda^3-(2a+e)\lambda^2+(a^2+b^2+2ae+\beta)%
\lambda-(a^2+b^2)e, 
\end{equation*}
where $\beta=-bc+kbdx_2^{k-1}$, we have the following result.

\begin{theorem}[stability regions for the symmetric period-two points]
\label{stregsym}\ Let the ACT family $F$ have the symmetric period-two
points $\pm p_2$, i.e., $k$ is odd and $bd(bc+c_{-1}(e))>0$, and let $J_{\text{sym}}(F_c)$ $($resp. 
$J_{\text{sym}}(F_{c,e})$$)$ denote
the stability region of the symmetric period-two orbit for $F_c$ $($resp. for $F_{c,e}$$)$. Then

\begin{enumerate}
\item If $bd<0$ then $J_{\text{sym}}(F_c)=J_{\text{sym}}(F_{c,e})=\o $. If $%
bd>0$, then $J_{\text{sym}}(F_c)\neq \o $ if and only if $J_{\text{tr}%
}(F_c)\neq \o $; moreover $J_{\text{non}}(F_{c,e})\neq \o $ if and only if $J_{\text{tr}}(F_{c,e})\neq 
\o $.

\item If $J_{\text{sym}}(F_{c,e})\neq \o $, then 
\begin{align*}
J_{\text{sym}}(F_{c,e})=&\Big\{(e,c)\in{\mathbb{R}^2}: \\
&\quad \frac{kc_{-1}(e)-\hat{c}(e)}{k-1}<-bc<\min\{\frac{kc_{-1}(e)+c_1(e)}{%
k-1}, c_{-1}(e)\}\Big\},
\end{align*}
where $\hat{c}(e)$, $c_{1}(e)$ and $c_{-1}(e)$ be in item $3$ of Theorem~\ref%
{stregorigin}.

\item For $F_{c,e}$, the point $(c_{-1}(e),e)$ is on the boundary of the
stability region of the symmetric period-two points  if and only if it is on the boundary of the stability
region of the origin.
\end{enumerate}
\end{theorem}

See Figure~3 for the regions of $J_{\text{sym}}(F_{c,e})$ together with $J_{%
\text{tr}}(F_{c,e})$ and $J_{\text{non}}(F_{c,e})$.

\begin{figure}[h]
\centerline{\hspace{0ex}\epsfig{file=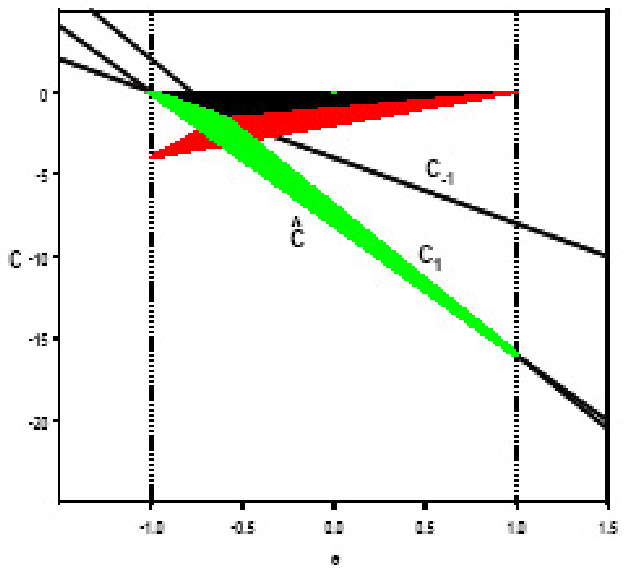,width=5.7cm,height=5.7cm}
\hspace{4ex}\epsfig{file=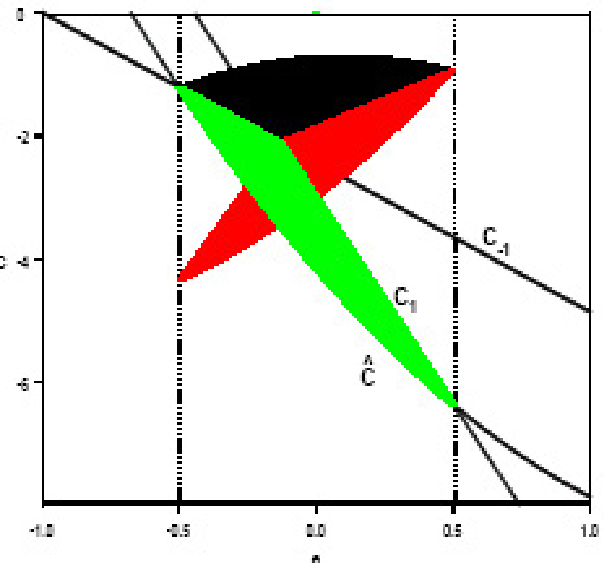,width=5.7cm,height=5.7cm}} \centerline{\ \
\ $(i)$\hspace{6.1cm}$(ii)$} 
\centerline{\hspace{0ex}\epsfig{file=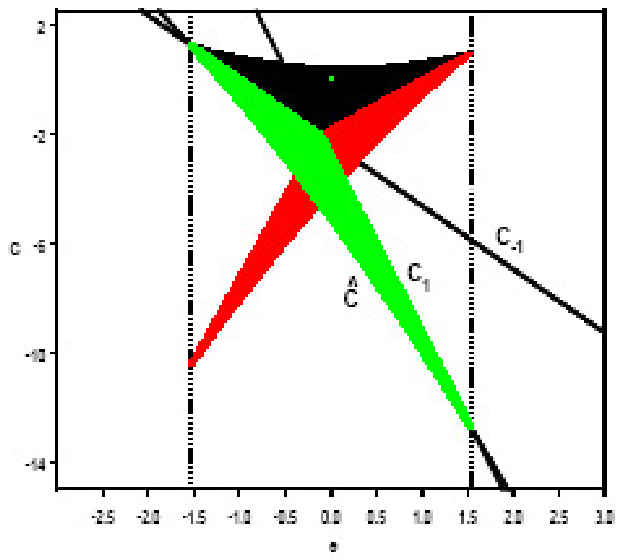,width=5.7cm,height=5.7cm}
} \centerline{\ \ \ $(iii)$}
\caption{The graphs of the bifurcation curves $c=\frac{kc_{-1}(e)-\hat{c}(e)%
}{-b(k-1)}$, $c=\frac{c_{-1}(e)}{-b}$ and $c=\frac{kc_{-1}(e)+c_{1}(e)}{%
-b(k-1)}$, indicated simply as $\hat{c}, c_{-1}$ and $c_{1}$, are shown in
the $(e,c)$-plane $("e"$ is the horizontal axis and $"c"$ is the vertical one%
$)$. The dashed lines in the figures are $e=\frac{\pm 1}{a^2+b^2}$. In
figures $1(i)$-$(iii)$, the stability regions $J_{\text{sym}}(F_{c,e})$ $($%
shaded in black$)$$, J_{\text{non}}(F_{c,e})$ $($shaded in gray$)$ and $J_{%
\text{tr}}(F_{c,e})$ $($shaded in light gray$)$ are shown for three cases:
when $a^2+b^2-1=0,>0$ and $<0$, all together with $b<0$, namely $(i)$ $a=0.6$
and $b=-0.8$, $(ii)$ $a=0.2$ and $b=-1.4$, and $(iii)$ $a=0.1$ and $b=-0.8$. 
}
\end{figure}

\begin{remark}
Note that for $k$ odd, the condition $\frac{bc-(1+e)((a+1)^2+b^2)}{bd}>0$
for existence of symmetric period-two points is satisfied precisely in the half-plane of $(e,c)$-plane which is
boundary by the straight line $c=c_{-1}(e)$ and lies to another side from this line than $J_{\text{tr}%
}(F_{c,e})$ the stability region of the trivial fixed point. So by using item $3$ of Theorem $7$ one map associate
appearance of these period-two points with period doubling bifurcation of the trivial fixed point, which occur at
the common segment $clo(J_{\text{tr}}(F_{c,e}))\cap clo(J_{\text{sym}}(F_{c,e}))$ of the line $%
c=c_{-1}(e)$ (recall that one of the eigenvalues of Jacobian matrix at the trivial fixed point is $-1$ as far as
parameter $(e,c)$ belong to this line).
\end{remark}

\section{Hyperbolic structure and infinitely many periodic points}

In~\cite{Du1985}, the first author of this paper proved that the set of
periodic points of the ACT family is contained in a bounded box. In \cite{LiMalkin2004}, the other two authors
extent the result as follows.

\begin{proposition}
\label{actnonwander} The nonwandering set of the ACT family $F$ lies in the
bounded box  
\begin{equation*}
\left\{(x,y,z)\in \mathbb{R^3}:\quad |x|\leq M,\quad |y|\leq \frac{%
a^2+b^2+2|a|+1}{|b|}M,\quad |z|\leq \frac{a^2+b^2+|a|}{|b|}M\right\}, 
\end{equation*}
where 
\begin{equation*}
M=\displaystyle \sqrt[k-1]{\frac{|a^2e+b^2e|+|a^2+b^2-bc+2ae|+|2a+e|+1}{|bd|}%
}.
\end{equation*}
\end{proposition}

\begin{remark}
\textrm{Note that if $f$ is a homeomorphism from $\mathbb{R^m}$ onto itself
then $f$ can be easily extended to  a homeomorphism $\bar f$ of the canonical one point compactification  $%
\bar{\mathbb{R}}^m:=\mathbb{R^m \bigcup \{\infty\}}$ with  $\bar
f(\infty)=\infty$. Therefore the topological entropy $h_{top}(f)$ can be
defined (see \cite{FriedlandMilnor1989})  as the usual topological entropy of a
homeomorphism of a compact set. On the other hand, using the results  of 
\cite{LiMalkin2004} on boundedness (and hence compactness) on $\Omega(f)$
under our assumptions, one may define  $h_{top}(f)$ to be $%
h_{top}(f|_{\Omega(f)})$. Since $h_{top}(\bar f)=h_{top}(\bar
f|_{\Omega(\bar f)})$ and  $\infty$ is an isolated nonwandering point for $%
\bar f$, the two definitions of the topological entropy of  $f$ agree. } 
\end{remark}

The following theorem shows the existence of horseshoe.

\begin{theorem}
Let $F_{a,e}$ be the ACT family with parameters $a$ and $e$. Suppose $b,c,d,k
$ satisfy

\begin{enumerate}
\item[$(i)$] $|b|>1$ and $|b-c|>1$
\end{enumerate}

and one of the following:

\begin{enumerate}
\item[$(ii1)$] $k$ is even, $\displaystyle bc>\frac{kb^2}{k-1}$,

\item[$(ii2)$] $k$ is even, $\displaystyle bc<\frac{kb^2}{k-1}(1-\sqrt[k-1]{k%
})$,

\item[$(ii3)$] $k$ is odd, $bd>0$, $\displaystyle bc>\frac{kb^2}{k-1}(1+\sqrt%
[k-1]{k})$,

\item[$(ii4)$] $k$ is odd, $bd<0$, $\displaystyle bc<\frac{kb^2}{k-1}(1-\sqrt%
[k-1]{k})$.
\end{enumerate}

Then for each pair of parameters $a$ and $e$ which are sufficiently close to 
$0$, $F_{a,e}$ has an invariant set, say $\Lambda_{a,e}$, such that the following properties hold:

\begin{enumerate}
\item $F_{a,e}$ has a hyperbolic structure on $\Lambda_{a,e}$.

\item $\Lambda_{a,e}$ is a Cantor set.

\item $F_{a,e}|\Lambda_{a,e}$ is topologically conjugate to the two-sided
shift on two $($resp. on three$)$ symbols for even $k$ $($resp. for odd $k$$)$.
\end{enumerate}
\end{theorem}

\begin{proof}By the persistence of hyperbolic invariant sets, we only need to consider the cease when $a=e=0$. We 
divide the ACT family into $8$ cases depending on the signs of $d$ and $bc$ for even $k$ and on the signs of $b$ and 
$d$ for odd $k$.  Since the proofs of all these $8$ cases are similar, we only give a proof of the case when $k>1$ 
is even, $d>0$, and $bc>0$.  Since 
\begin{align*}F_{0,0}: (x,y,z)\mapsto (\bar x, \bar y, \bar z)=(-b(y-z), bx, cx-dx^k),\end{align*} by eliminating 
$x$ among the expressions of $\bar y$ and $\bar z$, we obtains that the image of ${\Bbb R^3}$ under $F_{0,0}$ is the 
two-dimensional surface 
$$T=\left\{(x,y,z): x\in {\Bbb R}, z=\displaystyle \frac{c}{b}y-\frac{d}{b^k}y^k\right\}.$$  
By letting $\displaystyle\frac{dz}{dy}=0$ for $z=\displaystyle \frac{c}{b}y-\frac{d}{b^k}y^k$,  we have that the 
surface $T$ has the line 
$$\displaystyle \left\{(x, y_*, z_*): x\in {\Bbb R},  y_*=\sqrt[k-1]{\frac{b^{k-1}c}{kd}}, 
z_*=\frac{(k-1)c}{kb}\sqrt[k-1]{\frac{b^{k-1}c}{kd}}\right\}$$ as the set of extreme points.  By letting $z=0$, we 
get that $T$ intersects the $(x,y)$-plane at the lines $y=0$ and $y=r_*$ where $r_*=\displaystyle 
\sqrt[k-1]{\frac{b^{k-1}c}{d}}$.  It follows from the assumptions $d>0$ and $bc>0$ that $r_*>y_*>0$ and $z_*>0$.  
Also the assumption $bc>\displaystyle \frac{kb^2}{k-1}$ implies  $z_*>y_*$.  

Since $\bar{x}=-b(y-z)$, $F_{0,0}$ maps the infinite strip 
\begin{align*}S=\left\{(x,y,z): x\in {\Bbb R}, -r_*\leq z-y\leq 0\right\}\end{align*} onto the surface 
\begin{align*}F_{0,0}(S)=\left\{(x,y,z): z=\displaystyle \frac{c}{b}y-\frac{d}{b^k}y^k, \text{ $x$ lies between $0$ 
and $-br_*$}\right\},\end{align*} which is a subsurface of $T$.  Therefore, $S\cap F_{0,0}(S)$ consists of two 
disjoint connected sets whose projection to the $(y,z)$-plane are shown in (i) of Figure~4; (ii) of Figure~4 is for 
the case when $k$ is odd.  

\begin{figure}[h]
\vspace{0.3cm}
\centerline{\epsfig{file=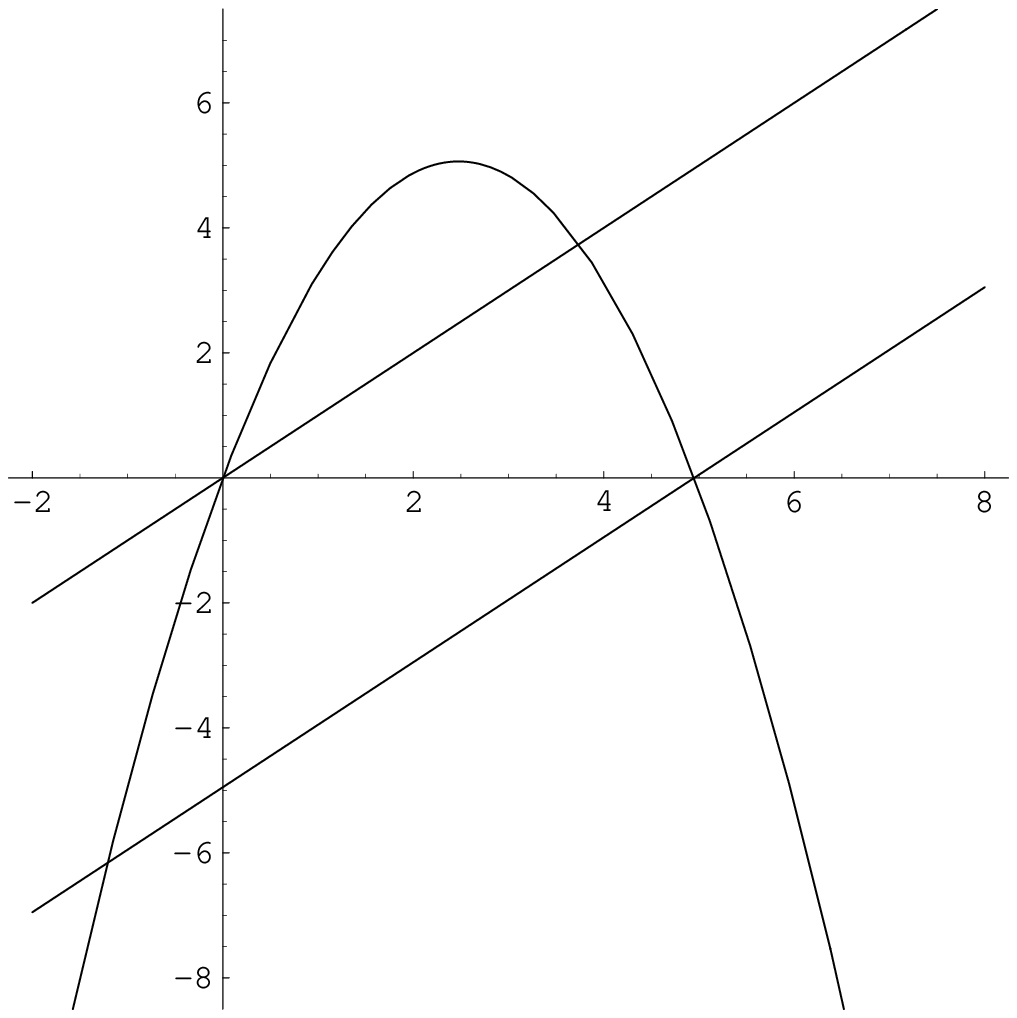,width=3.5cm,height=5cm}
\hspace{6ex}\epsfig{file=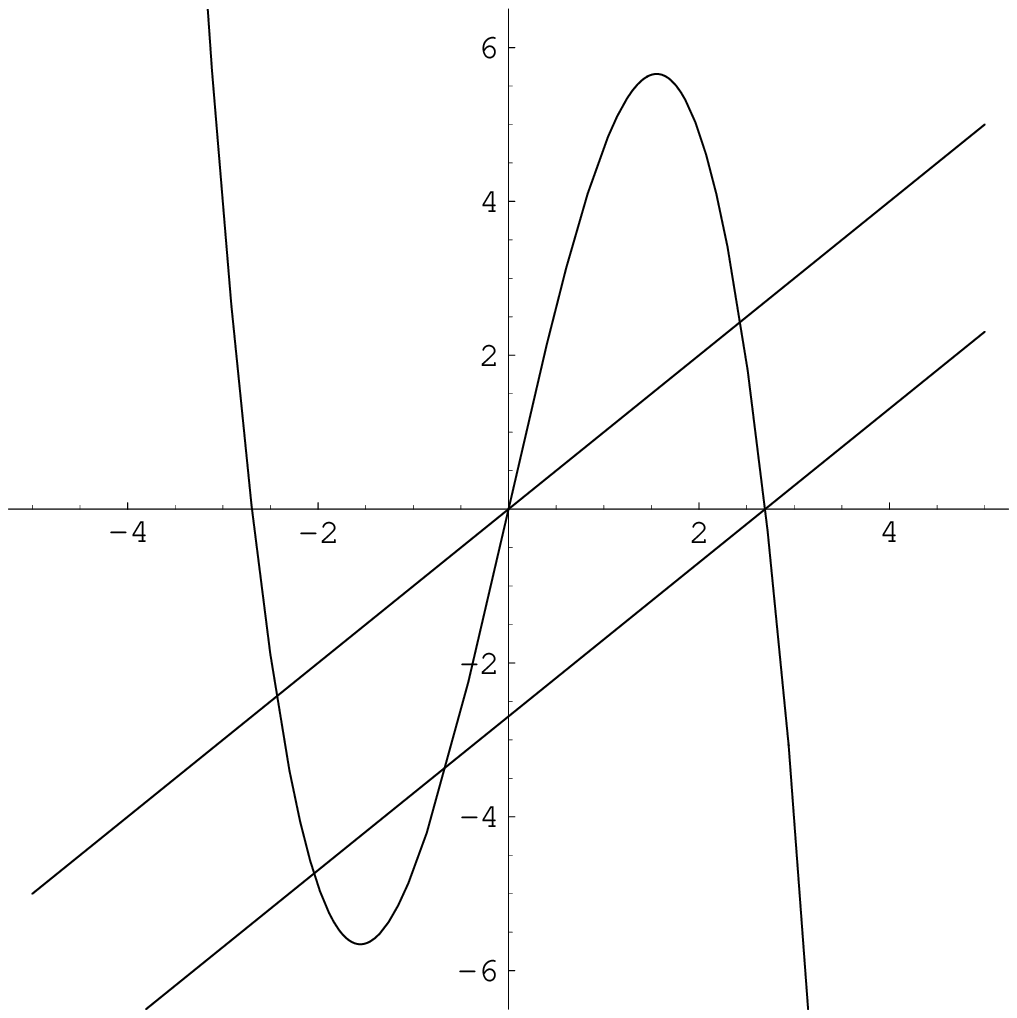,width=3.5cm,height=5cm}
}
\centerline{\ \ \ $(i)$\hspace{4.1cm}$(ii)$}
\caption{The boundaries of $S$ and $F(S)$.}
\end{figure}

Next, we study $F_{0,0}^{-1}(S)\cap S$, where $F_{0,0}^{-1}$ denotes the preimage of $S$ under $F_{0,0}$.  
Considering $\bar{z}-\bar{y}=-r_*$ and  $\bar{z}-\bar{y}=0$, one has the equations $d x^k+(b-c)x-r_*=0$ and $d 
x^k+(b-c)x=0$.  The two equations give four distinct real roots in $x$, namely  
$\lambda_1<\lambda_2<\lambda_3<\lambda_4$.  Then  
\begin{align*}
F_{0,0}^{-1}(S)&=\{(x,y,z): \lambda_1\leq x\leq \lambda_2, y\in {\Bbb R}, z\in 
{\Bbb R}\}\\
&\quad \bigcup \{(x,y,z): \lambda_3\leq x\leq \lambda_4, y\in {\Bbb R}, z\in {\Bbb R}\}.
\end{align*} 

Let $\Lambda_0=\displaystyle \bigcap_{n=-\infty}^\infty F_{0,0}^n(S)$.  We need to check the hyperbolicity 
conditions for points in a small neighborhood $U$ of $F_{0,0}^{-1}(S)\cap S$. First of all, we make a change of 
coordinate via $u=y-z$ and $v=y+z$.  Then in the new coordinate $(x,u,v)$, we have 
\begin{align*}F_{0,0}(x,u,v)=(-bu, bx-cx+dx^k, bx+cx-dx^k).
\end{align*}
Let ${\bf x}=(x,y)\in {\Bbb R^2}$ and write $F_{0,0}(x,u,v)=(G({\bf x},{v}), H({\bf x},{v}))$, where $G:{\Bbb 
R^3}\to{\Bbb R^2}$ and $H:{\Bbb R^3}\to{\Bbb R}$.  Then the sufficient conditions of hyperbolicity are the following 
(refer to  \cite{AfraimovichBykovShilnikov1983}): 
\begin{align}
&\|\big(\frac{\partial G}{\partial {\bf x}}\big)^{-1}\|<1,\\
&\|\frac{\partial H}{\partial v}\|<1,\\
&1-\|\big(\frac{\partial G}{\partial {\bf x}}\big)^{-1}\|\ \cdot \|\frac{\partial H}
{\partial v}\|>2\sqrt{\|\frac{\partial H}{\partial {\bf x}}\cdot 
\big(\frac{\partial G}{\partial {\bf x}}\big)^{-1}\|\ \cdot  \|G_{ v}\|\ \cdot \|\big(\frac{\partial G}{\partial 
{\bf x}}\big)^{-1}\|},\\
&(1-\|\big(\frac{\partial G}{\partial {\bf x}}\big)^{-1}\|)\cdot (1-\|\frac{\partial H}{\partial 
v}\|)>\|\frac{\partial H}{\partial {\bf x}}\cdot 
\big(\frac{\partial G}{\partial {\bf x}}\big)^{-1}\|\ \cdot \|\frac{\partial G}{\partial v}\|,
\end{align}
where $\|\cdot\|=\displaystyle\sup_{{(\bf x,v)}\in U}|\cdot|$.

The Jacobian matrix of $F_{0,0}$ in the new coordinate is 
$$\begin{bmatrix}\displaystyle\frac{\partial G}{\partial {\bf x}}& \displaystyle\frac{\partial G}{\partial v}\\ 
\ &\ \\
\displaystyle\frac{\partial H}{\partial {\bf x}}&\displaystyle\frac{\partial H}{\partial v} \end{bmatrix}=
\begin{bmatrix}0&-b&0\\ b-c+kdx^{k-1}&0&0 \\
\\ b+c-kdx^{k-1}&0&0\end{bmatrix}.$$
Since $\|\frac{\partial H}{\partial v}\|=\|\frac{\partial G}{\partial v}\|=0$, $(7)$ and $(8)$ are clear and $(9)$ 
follows from $(6)$.  We only need to check $(6)$.  By the definition of $\lambda_i$'s and the derivative of 
$f(x)=dx^k+(b-c)x$, we have that  $$\min\{|b-c+kdx^{k-1}|: x\in 
[\lambda_1,\lambda_2]\cup[\lambda_3,\lambda_4]\}=|b-c|.$$ 
Thus (6) follows from the assumptions $|b|>1$ and $|b-c|>1$.  
\end{proof}

\section{Chaos near an anti-integrable limit}

In~\cite{LiMalkin2004diff}, the second and third authors consider solutions
of families of difference equations $\Phi_\lambda(y_n, y_{n+1},\dots,
y_{n+m})=0$, $n\in \mathbb{Z}$, with $m+1$ real variables for parameters $%
\lambda$ near those exceptional values $\lambda_0$ for which the difference
equation function depends on only one variable: $\Phi_{\lambda_0}(y_0,%
\dots,y_m)=\varphi(y_N)$  with some integer $0\leq N\leq m$ and function $%
\varphi$. It is proved there that if $\varphi$ has $k$ simple zeros, then
among  solutions for parameter values close to exceptional ones , there are
"topological $k$-horseshoes", i.e., solutions for which the restriction of
the shift map  is conjugate to the two-sided full shift on $k$ symbols, and
if the difference equations correspond to smooth diffeomorphisms on $\mathbb{%
R^m}$, then the topological $k$-horseshoes persist among orbits of such
diffeomorphisms. For the exceptional parameter values, no real maps are
correspondent, and such a situation reminds the so called anti-integrable
limit approach inspired by Aubry and Abramovici's anti-integrable concept in~%
\cite{AubryAbramovici1990}; see also~\cite{ChenYC2004, MackayMeiss1992,
Qin2001}.

For an initial point $p=(x_0,y_0,z_0)$, denote the $n$-th iteration of $p$
under the ACT map $F$ by $(x_n,y_n,z_n)$. Then for any $n\in{\mathbb{Z}}$, we have the following system consisting 
of $%
7$ equations  (the index $i$ for labeling the equations bellow takes the
values $1$ and $2$): 
\begin{align*}
\begin{cases}
x_{n+i}=ax_{n+i-1}-b(y_{n+i-1}-z_{n+i-1}), & \hspace{6cm}(1,i) \\ 
y_{n+i}=bx_{n+i-1}+a(y_{n+i-1}-z_{n+i-1}), & \hspace{6cm}(2,i) \\ 
z_{n+i}=cx_{n+i-1}-dx_{n+i-1}^k+ez_{n+i-1}, & \hspace{6cm}(3,i) \\ 
x_{n+3}=ax_{n+2}-b(y_{n+2}-z_{n+2}). & \hspace{6cm}(1,3) \\ 
& 
\end{cases}%
\end{align*}
From equations~$(1,1)$ and $(2,1)$, we can express $y_{n+1}$ in terms of $x_n
$ and $x_{n+1}$. Then plugging $y_{n+1}$ into equation~$(1,2)$, we can express $z_{n+1}$ in terms of $x_{n}$%
, $x_{n+1}$ and $x_{n+2}$. The expressions are 
\begin{align*}
y_{n+1}&=\frac{a^2+b^2}{b}x_{n}-\frac{a}{b}x_{n+1}, \\
z_{n+1}&=\frac{a^2+b^2}{b}x_{n}-\frac{2a}{b}x_{n+1} +\frac{1}{b}x_{n+2}.
\end{align*}
Similarly, from equations~$(1,2)$, $(2,2)$ and $(1,3)$, we can express $%
z_{n+2}$ in terms of $x_{n+1}$, $x_{n+2}$ and $x_{n+3}$. Then, plugging the expressions of $z_{n+1}$ and $z_{n+2}$
into equation~$(3,2)$, we get the difference equation 
\begin{align*}
dx_{n+1}^k+\displaystyle -\frac{a^2e+b^2e}{b}x_n+\frac{a^2+b^2-bc+2ae}{b}%
x_{n+1}-\frac{2a+e}{b}x_{n+2} +\frac {1}{b}x_{n+3}=0.
\end{align*}

Suppose the parameters $a,b$ and $e$ to be fixed and let $c\to\infty$, $d\to
\infty$ in such a way that $\frac{d}{c}=\text{constant}:=A>0$. Denote $\lambda=\frac{1}{c}$, then $d=%
\frac{A}{\lambda}$ and for the $x$-coordinate of orbit under the ACT map $f=f_\lambda$, we have the
difference equation 
\begin{align}  \label{actdif}
x_{n+1}(x_{n+1}^{k-1}-A) +\lambda\left(\displaystyle -\frac{a^2e+b^2e}{b}x_n+%
\frac{a^2+b^2+2ae}{b}x_{n+1}-\frac{2a+e}{b}x_{n+2}+\frac {1}{b}%
x_{n+3}\right) =0.
\end{align}
The difference equation~\eqref{actdif} corresponds to the map $f_\lambda$
for $\lambda\neq 0$ because between solutions $\underline{x}=(x_n)_{n=-\infty}^\infty$ of~\eqref{actdif} and
full orbits $\underline{p}=(p_n)_{n=-\infty}^\infty$ of $f_\lambda$, we have a conjugacy 
$\underline{x}\mapsto \underline{p}$, given by 
\begin{align*}
p_n=(x_n,\frac{a^2+b^2}{b}x_{n-1}-\frac{a}{b}x_{n}, \frac{a^2+b^2}{b}x_{n-1}-%
\frac{2a}{b}x_{n} +\frac{1}{b}x_{n+2}),
\end{align*}
while the inverse $\underline{p}\mapsto \underline{x}$ is given by $%
x_n=\pi_1(p_n)$. For the limit value $\lambda_0=0$ of parameter (in this case, it is an anti-integrable limit;
refer to~\cite{AubryAbramovici1990}), we have 
\begin{align*}
\psi(x_{n+1}) :=x_{n+1}(x_{n+1}^{k-1}-A).
\end{align*}
Thus the function $\psi$ has at least two simple zeros; more precisely $\psi$
has two simple zeros $\{0,\sqrt[k-1]{A}\}$ when $k$ is even, and three simple zeros $\{0,\pm\sqrt[%
k-1]{A}\}$ when $k$ is odd. By applying Theorem 3 of~\cite{LiMalkin2004diff}%
, we have the following results.

\begin{proposition}
Suppose that for the family of the ACT maps~\eqref{actmap}, the parameters $%
a,b$, and $e$ are fixed while $c\to \infty$ and $d\to\infty$ in such a way that $\frac{d}{c}=\text{constant}>0$.
Let $\lambda=\frac{1}{c}$. Then for all $|\lambda|$ sufficiently small, the ACT map $F=F_\lambda$ has a closed
invariant set $\Lambda_\lambda$ such that $F_\lambda|\Lambda_\lambda$ is conjugate to the full shift on either two 
or
three symbols depending on whether $k$ is even or odd, respectively.
\end{proposition}

There are also other parameter routes in families of Arneodo-Coullet-Tresser
maps, for which the above arguments apply. Namely, if we consider $a$ and $c$
to be fixed while $e\to 0$, and $b\to \infty$ and $d\to\infty$ in such a way
that $\frac{b}{d}=\text{constant}>0$, then similarly to the lines of the
previous proposition, we will have the following.

\begin{proposition}
Suppose that for a family of the ACT maps~\eqref{actmap}, the parameters $a$
and $c$ are fixed while $e\to 0$, and $b\to \infty$ and $d\to\infty$ in such a way that $\frac{b}{d}=\text{constant%
}>0$. Let $\lambda=\frac{1}{d}$.  Then for all sufficiently small $|e|$ and $%
|\lambda|$, the ACT map $F=F_{e,\lambda}$ has a closed invariant set  $%
\Lambda_{e,\lambda}$ such that $F_{e,\lambda}|\Lambda_{e,\lambda}$ is
conjugate to the full shift on either two or three  symbols depending on whether $k$ is even or odd, respectively.
\end{proposition}

\section{Numerical results}

Fix $a,b,d,k$ such that the stability regions of the origin, the nontrivial
fixed point $p_1$ and the symmetric period-$2$ point $p_2$ are not empty. For any $e$ we denote the boundaries
of the stability regions as follows: $\hat{c}_{\text{tr}}(e)=\frac{\hat{c}(e)%
}{-b}$, $c_{-1,\text{tr}}(e)=\frac{c_{-1}(e)}{-b}$, $\hat{c}_{\text{sym}}(e)=%
\frac{kc_{-1}(e)-\hat{c}(e)}{-b(k-1)}$, $c_{-1,\text{non}}(e)=\frac{kc_1(e)+c_{-1}(e)}{b(k-1)}$, and 
$\hat{c}_{\text{%
non}}(e)=\frac{k c_1(e)+\hat{c}(e)}{b (k-1)}$, where $\hat{c}(e)$, $c_1(e)$
and $c_{-1}(e)$ are in item $3$ of Theorem $5$; see Figure 5. Let $e_{\pm}=\pm\frac{1}{a^2+b^2}$, then 
$\hat{c}_{\text{tr}}(e_{\pm})=\frac{%
(a^2+b^2-1)[(ae_\pm-1)^2+b^2e_\pm^2]} {b}$ and the Jacobian matrix of $%
F_{c_1(e_\pm),e_\pm}$ at all fixed points has all eigenvalues lying on the
unit circle. For $e$ with $e_-<e<e_+$, the Jacobian matrix of $F_{\hat{c}_{\text{%
tr}}(e), e}$ (resp. $F_{\hat{c}_{\text{sym}}(e), e}$ and $F_{\hat{c}_{\text{non}}(e),e}$) at the
trivial fixed point (resp. at the symmetric period-$2$ point and at the nontrivial fixed points) has simple
complex eigenvalues lying on the unit circle and the Jacobian matrix of $F_{c_{-1,\text{tr}}(e), e}$ (resp. $%
F_{c_{-1,\text{non}}(e), e}$) at the trivial (resp. nontrivial) fixed points has simple real eigenvalue at $-1$.
Numerical results indicate the following:

\begin{enumerate}
\item The family $F_{c,e}$ undergoes Hopf bifurcations at $(c,e)$ equal to $(%
\hat{c}_{\text{tr}}(e), e)$, $(\hat{c}_{\text{sym}}(e),e)$ and $(\hat{c}_{\text{non}}(e),e)$ for almost
all $e$ with $e_-<e<e_+$. As the value of $c$ decreases, the former bifurcation is subcritical and the invariant
circle (before disappeared) is saddle type while the latter two are supercritical and the appeared invariant circle 
is
asymptotically stable. It can be shown that resonance cases take places for several parameters; e.g., for $c$ equal
to $\hat{c}_{\text{tr}}(e)$ and $\hat{c}_{\text{non}}(e)$, $F_{c,e}$ has resonance $(1:4)$ at $e=\frac{2a}{%
a^2+b^2-1}$ and resonance $(1:3)$ at $e=\frac{2a+1}{a^2+b^2-1}$. See Figure 5.

\item The family $F_{c,e}$ undergoes period-doubling bifurcations at $(c,e)$
equal to $(c_{-1,\text{tr}}(e),e)$ and $(c_{-1,\text{non}}(e), e)$ for all $e$ with $e_-<e<\bar{e}$, where $\bar{e}$
with $c_{-1,\text{tr}}(\bar{e})=c_{-1, \text{non}}(\bar{e})$. In fact, one
can show that the symmetric period-two points $\pm p_2$ bifurcates from the 
origin at $(c_{-1,\text{tr}}(e),e)$. See Figure 5.

\item If $e=e_\pm$, then for each $c$ near $\hat{c}_{\text{tr}}(e_\pm)$ with 
$d(c-\hat{c}_{\text{tr}}(e_\pm))>0$, the map $F_{c}$ has infinitely many invariant tori (or double tori); one
lies inside another and none of these tori is asymptotically stable. See Figure $6(i)$ with $a=0.2, b=-1.4, c=-0.94,
d=-1, e=0.5$ and $k=3$ and Figure $6(ii)$ with $a=0.6, b=0.8, c=-0.01, d=-1, e=-1$ and $k=4$. For the former case, 
one
invariant circle is drawn and it indicates that all other invariant tori are bifurcated from the invariant
circle through a degenerate Hopf bifurcation by taking a Poincare section. For the latter case, each double
torus consists of two tori which are mapped to each other alternatively and is invariant under the map $F_c$.

\item For some parameters, e.g., $a=b=0.5$, $d=e=1$, and $k=3$, the family $%
F_c$ undergoes sequences of wind-doubling bifurcations of invariant circle, and has coexistence of
stable invariant circles and strange attractors for some parameters $c$. See Figure 7.

\item For some parameters, e.g., $a=e=0.01$, $b=1.1$, $c=3.6578$, $d=1$, and 
$k=3$, the family $F$ has a strange attractor. See Figure 8.
\end{enumerate}

\begin{figure}[h]
\centerline{\epsfig{file=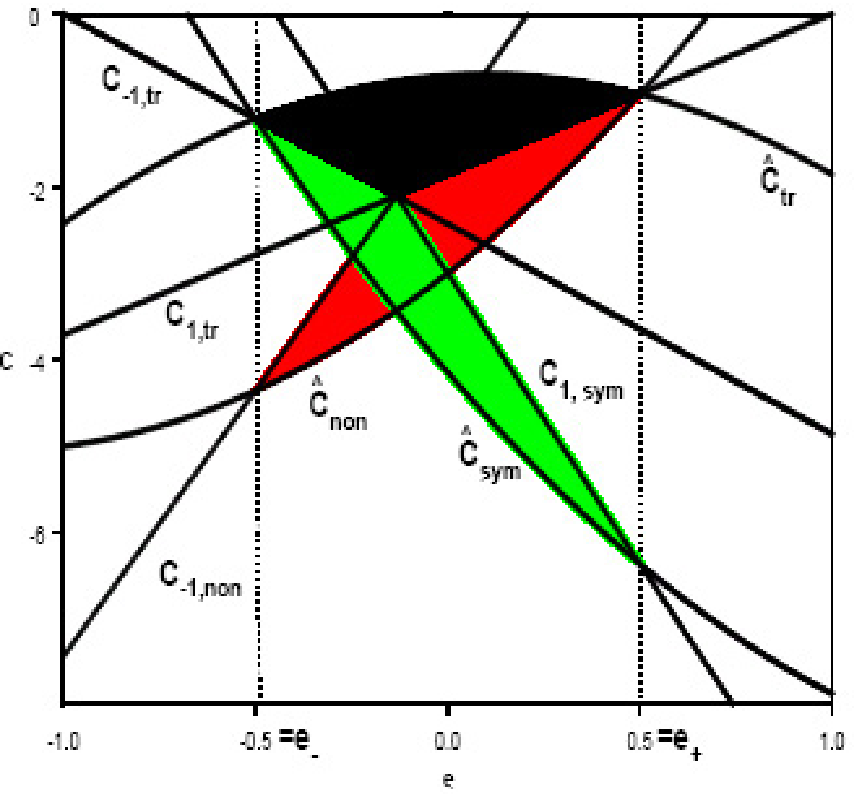,width=6cm,height=6cm}
}
\caption{The bifurcation curves.}
\end{figure}

\begin{figure}[h]
\centerline{\epsfig{file=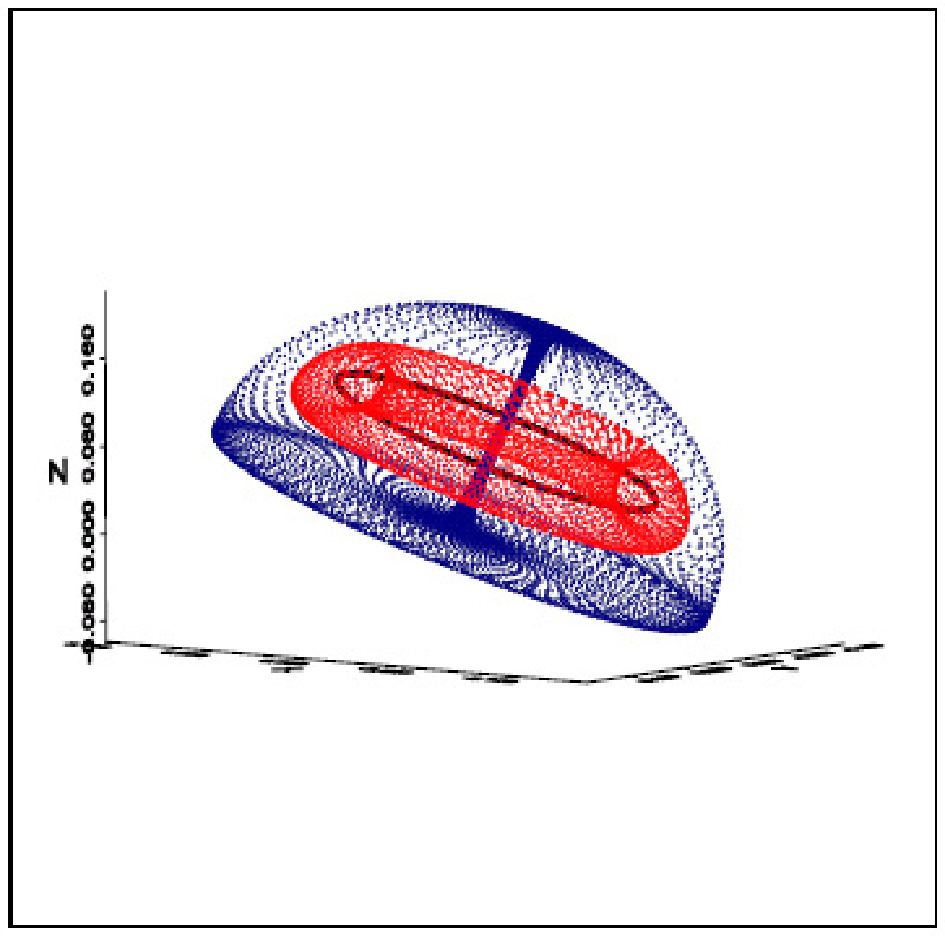,width=6cm,height=6cm}
\hspace{2ex}\epsfig{file=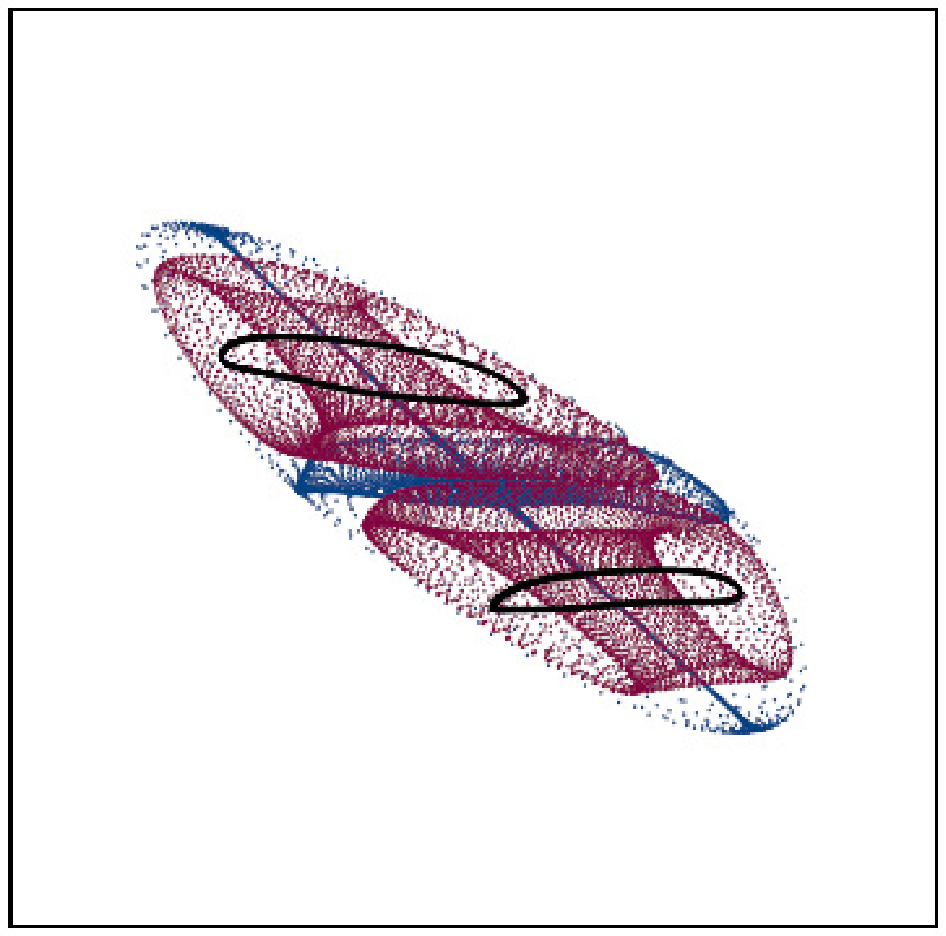,width=6cm,height=6cm}
} \centerline{$(i)$\hspace{6.2cm}$(ii)$}
\caption{Invariant tori: one lies inside another for $(i)$ $k=3$ and $(ii)$ $%
k=4$.}
\end{figure}

\begin{figure}[h]
\centerline{\epsfig{file=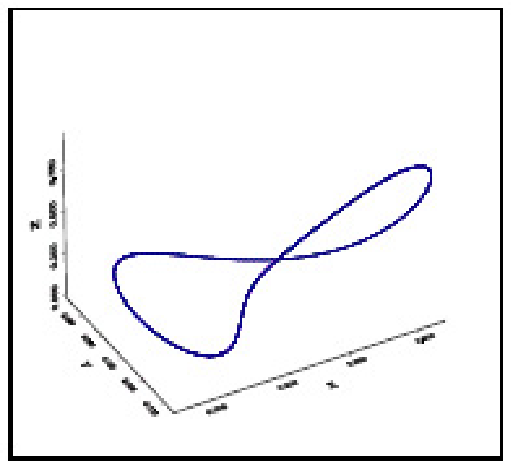,width=4cm,height=4cm}
\hspace{2ex}\epsfig{file=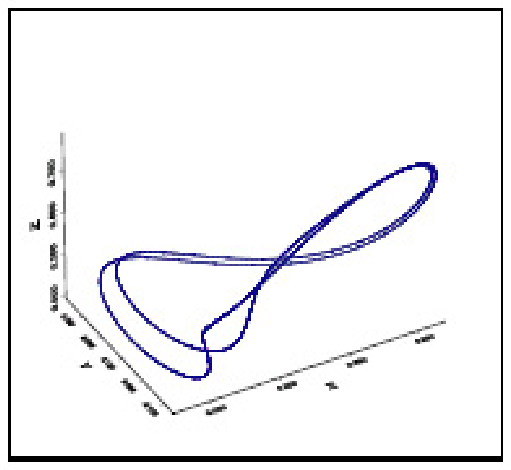,width=4cm,height=4cm}
\hspace{2ex}\epsfig{file=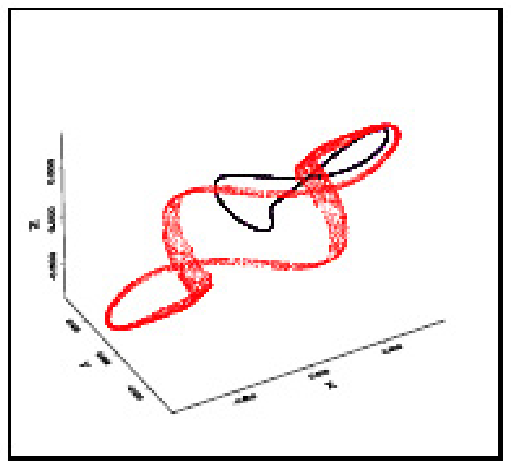,width=4cm,height=4cm}
} \centerline{$(i)$\hspace{3.7cm}$(ii)$\hspace{3.7cm}$(iii)$}
\caption{$(i)$ and $(ii)$ wind-doubling of invariant circle, and $(iii)$
coexistence of invariant circle and strange attractor.}
\end{figure}

\begin{figure}[h]
\centerline{\epsfig{file=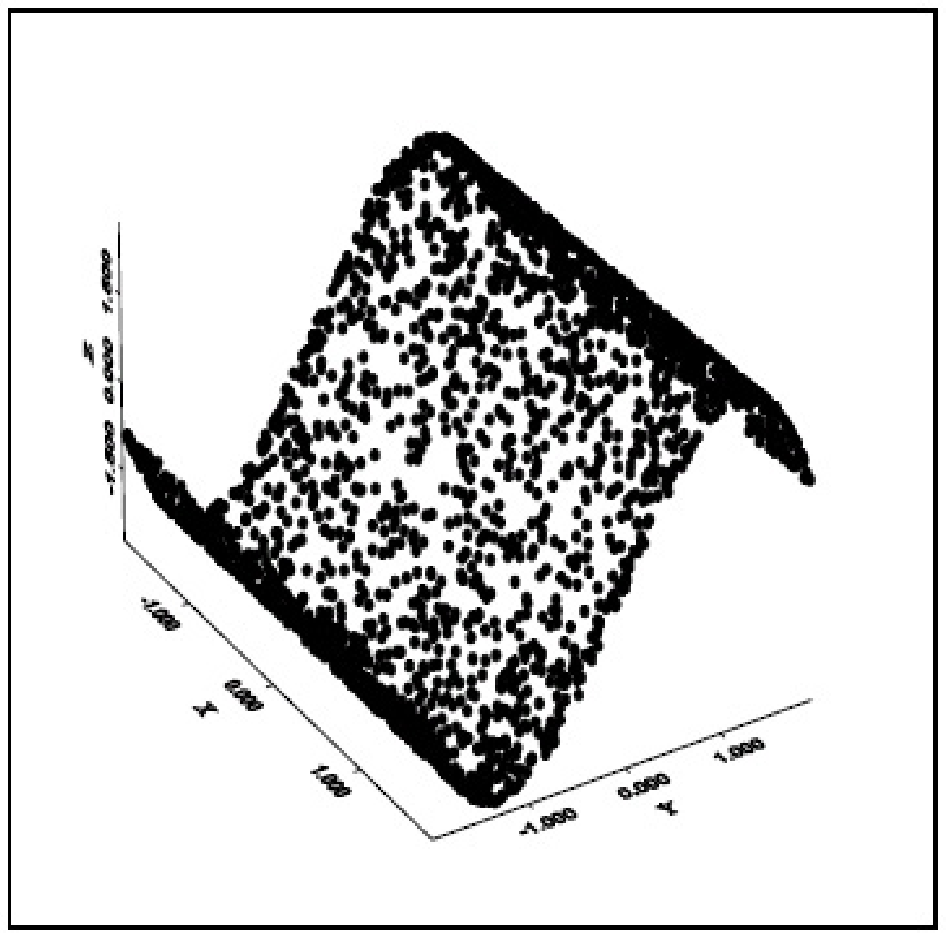,width=8cm,height=8cm}}
\caption{Existence of strange attractor.}
\end{figure}

\bigskip \noindent \textbf{{\Large Appendix: Proofs of Propositions~\ref%
{polystable} and ~\ref{bifstable}}}\ \newline

\begin{proof}[Proof of Proposition~\ref{polystable}]
First, we consider the case when $D=0$.  Let $Q(\lambda)=\lambda^2+A\lambda+B$, then $P(\lambda)=\lambda 
Q(\lambda)$.   For the "only if" part of the proposition, the stability of $P$ implies that the product of two roots 
of $Q$ has absolute value less than one and so $B<1$, i.e., $\hat{\alpha}>0$.  Together with the fact that the 
leading coefficient of $Q$ is positive, we have that $P(1)=Q(1)>0$ and $-P(-1)=Q(-1)>0$.  For the "if" part, by the 
assumptions, we have that $B<1$, $1+A+B=Q(1)=P(1)>0$, and $1-A+B=Q(-1)=-P(-1)>0$.  Adding the last two inequalities 
gives that $B>-1$ and so together with the first inequality we obtain $|B|<1$, that is, the product of two roots of 
$Q$ has absolute value less than one.  It follows that the two roots of $Q$ lie inside the unit circle because of 
the following facts: $Q(1)=P(1)>0$, $Q(-1)=-P(-1)>0$, and that the leading coefficient of $Q$ is positive.  The 
proof of the proposition in the case when $D=0$ is completed.   

Next, we consider the case when $D\neq 0$.  Before the proof, we give some preparations.  Define 
$P_\alpha(\lambda)=P(\lambda)+\alpha\lambda$ as a one-parameter family of polynomials with parameter $\alpha\in 
{\Bbb R}$.  
Then for any fixed $\lambda_0\neq 0$, $P_\alpha(\lambda_0)$ has the {\em monotonicity property} in $\alpha$: 
if $\lambda_0>0$ ($\lambda_0<0$, $\lambda=0$ respectively), then  $P_\alpha (\lambda_0)$ increases 
(decreases, equals to $D$, respectively) as $\alpha$ increases.
 It is clear that there is a unique $\alpha$, namely $\alpha=-P(1)$, such that $P_\alpha(1)=0$. Similarly, 
$P_\alpha(-1)=0$ for a unique $\alpha$, namely $\alpha=P(-1)$.  Since the product of the three roots of $P_\alpha$ 
is $-D$, it follows that there is a unique $\alpha$, namely $\alpha=\frac{P(-D)}{D}=\hat{\alpha}$, such that 
$P_\alpha$ has two roots whose product is equal to one.  In particular, if $P_\alpha$ has two complex conjugate 
roots on the unit circle, then $\alpha=\hat{\alpha}$.  For convenience, we denote $\alpha_1=-P(1)$, 
$\alpha_{-1}=P(-1)$ 
and  $\alpha_{max}=\max\{\alpha_1,\alpha_{-1},\hat{\alpha}\}$.
Note that the numbers $\alpha_1,\alpha_{-1}, \hat{\alpha}$ partition the real line into disjoint open intervals and 
for any $\alpha,\alpha^\prime$ in the same interval of this partition, we have that $P_\alpha$ and 
$P_{\alpha^\prime}$ are either stable or unstable simultaneously.  Indeed, this follows from the fact that roots of 
polynomials depends continuously on the coefficients.  Because of this, we call $\alpha_1$, $\alpha_{-1}$ and 
$\hat{\alpha}$ the {\em bifurcation values} of the family $P_\alpha$.

Now, we prove the "only if" part.  The first statement $|D|<1$ follows immediately from the stability of $P$.  To 
prove the second statement, we need the following two claims:

Claim (i) {\em $P_\alpha$ is unstable for all $\alpha\geq \alpha_{max}$. } 
Suppose the claim is not true. Then there is $\alpha^\prime>\alpha_{max}$ such that $P_{\alpha^\prime}$ is stable.  
Thus for each $\alpha>\alpha_{max}$ the polynomial $P_\alpha$ is stable and so $\lambda_\alpha:=\min\{\lambda\in 
(-1,1): \lambda$ is a root of $P_\alpha\}$ is well defined.  Thus 
\begin{align*}
|\alpha\lambda_\alpha|=|-P(\lambda_\alpha)|
=|\lambda_\alpha^3+A\lambda_\alpha^2+B\lambda_\alpha+D|\leq 1+|A|+|B|+|D|.
\end{align*}  It follows that $\lambda_\alpha\to 0$ as $\alpha\to+\infty$.  Thus for sufficiently large $\alpha$, 
the product of the two roots of $P_\alpha$ other than $\lambda_\alpha$ is $\frac{-D}{\lambda_\alpha}$, whose 
absolute value is large.  
This contradicts the stability of $P_\alpha$.  

Claim (ii) {\em $P_\alpha$ is unstable for all $\alpha\leq \max\{\alpha_1, \alpha_{-1}\}$}.  Indeed, by taking 
$\lambda_0=1$, the monotonicity property in $\alpha$ of $P_\alpha(\lambda_0)$ implies that 
$P_\alpha(1)<P_{\alpha_1}(1)=0$ for all $\alpha<\alpha_1$.  On the other hand, the positivity of the leading 
coefficient of $P_\alpha$ implies that  $P_\alpha(\lambda)>0$ for some $\lambda$ large.  Thus for all 
$\alpha<\alpha_1$ the polynomial $P_\alpha$ has a real root bigger than $1$ and so $P_\alpha$ is unstable.  The 
proof of the fact that $P_\alpha$ is unstable for all $\alpha< \alpha_{-1}$ is similar by taking $\lambda_0=-1$.  

Combining the claims (i) and (ii), we get that if $\hat{\alpha}\leq \max\{\alpha_1,\alpha_{-1}\}$ then $P_\alpha$ is 
unstable for all $\alpha\in {\Bbb R}$.  In other words, if $P_\alpha$ is stable then $\max\{\alpha_1, 
\alpha_{-1}\}<\alpha<\hat{\alpha}$. In particular, by the stability of $P_\alpha$ with $\alpha=0$, we have that 
$\max\{\alpha_1, \alpha_{-1}\}<0<\hat{\alpha}$.  

Next, we prove the "if" part.  Since the only bifurcation values of $P_\alpha$ with parameter $\alpha$ are 
$\alpha_1, \alpha_{-1}$ and $\hat{\alpha}$, it is sufficient to show that $P_\alpha$ is stable for some 
$\alpha^\prime$ with $\max\{\alpha_1, \alpha_{-1}\}<\alpha^\prime<\hat{\alpha}$.  By the assumption $|D|<1$ and the 
definition of $\hat{\alpha}$, for $\alpha=\hat{\alpha}$ the polynomial $P_{\alpha}$ has a real root at $-D$ which 
lies in $(-1,1)$.  Since $\hat{\alpha}>\max\{\alpha_1,\alpha_{-1}\}$ and the monotonicity property of 
$P_\alpha(\lambda_0)$ with $\lambda_0=1$ and $ -1$ imply $P_{\hat{\alpha}}(-1)<0<P_{\hat{\alpha}}(1)$.  Now take an 
$\alpha^\prime$ slightly less than $\hat{\alpha}$ so that $\max\{\alpha_1, \alpha_{-1}\}<\alpha^\prime<\hat{\alpha}$ 
and $P_{\alpha^\prime}(-1)<0<P_{\alpha^\prime}(1)$. By the monotonicity property of $P_\alpha(\lambda_0)$ with 
$\lambda_0=-D$, we get that $P_{\alpha^\prime}$ has a real root, say $\lambda^\prime$, with 
$|-D|<|\lambda^\prime|<1$.  Let $\lambda_1,\lambda_2$ denote the other two roots of $P_\alpha$.  Since 
$\lambda^\prime\lambda_1\lambda_2=-D$, we have $|\lambda_1\lambda_2|<1$.  So if $\lambda_1$ and $\lambda_2$ are a 
pair of complex conjugate, then $P_{\alpha^\prime}$ is stable.  If $\lambda_1,\lambda_2$ are real numbers, by using 
the inequalities $|\lambda_1\lambda_2|<1$ and $P_{\alpha^\prime}(-1)<0<P_{\alpha^\prime}(1)$ and the fact that 
$P_{\alpha^\prime}$ is a cubic polynomial with positive leading coefficient, we have that $\lambda_1,\lambda_2$ lie 
in the interval $(-1,1)$ and so $P_{\alpha^\prime}$ is stable. 
\end{proof}

\begin{proof}[Proof of Proposition~\ref{bifstable}]
From the proof of Proposition~\ref{polystable}, we have known that $-P_0(1), P_0(-1)$ and 
$\hat{\alpha}_0:=-D^2+AD-B+1$ are the only bifurcation values of the family $P_\beta$.
Applying Proposition~\ref{polystable} to the polynomial, we get that  $P_\beta$ is stable if and only if  
\begin{align}
\text{$|D|<1$ \hspace{3ex} and \hspace{3ex} $\max\{-P_\beta(1), P_\beta(-1)\}<0<\hat{\alpha}(P_\beta)$.} \end{align}
It is easy to see that $-P_\beta(1)=-P_0(1)-\beta$, $P_\beta(-1)=P_0(-1)-\beta$ and 
$\hat{\alpha}(P_\beta)=\hat{\alpha}_0-\beta$.  Hence, (1) implies that $P_\beta$ is stable if and only if 
\begin{align}
\text{$|D|<1$ \hspace{3ex} and \hspace{3ex} $\max\{-P_0(1), P_0(-1)\}<\beta<\hat{\alpha}_0$.} 
\end{align}

From (2), we have that the stable interval of $P_\beta$ exists if and only if
\begin{align}
\text{$|D|<1$, \hspace{2ex} $\hat{\alpha}_0+P_0(1)>0$ \hspace{2ex} and \hspace{2ex} $\hat{\alpha}_0-P_0(-1)>0$.} 
\end{align}
By simple computations, we get 
\begin{align*}
\hat{\alpha}_0+P_0(1)&=(D+1)(2+A-D)\text{ and }\\
\hat{\alpha}_0-{P_0(-1)}&=(1-D)(2-A+D).
\end{align*}
Therefore, (3) is equivalent to 
\begin{align*}
\text{$|D|<1$ \hspace{2ex}  and \hspace{2ex} $|A-D|<2$.} 
\end{align*}
The equivalence of $(a)$ and $(b)$ is completed.  The equivalence of $(b)$ and $(c)$ was actually established in the 
proof of Proposition~\ref{polystable}. 

The first statement in item $2$ follows immediately from item $1$ and conditions (2).  Finally, we prove the second 
statement in item $2$.
By the definition of $\hat{\alpha}_0$, the product of two roots of $P_{\hat{\alpha}_0}$, say $\lambda_1$ and 
$\lambda_2$, is equal to one.  Since $P_\beta$ is stable for $\beta$ slightly less than $\hat{\alpha}_0$, we have 
that $|\lambda_1|=|\lambda_2|=1$.  Because $\lambda_1$ and $\lambda_2$ are different from $1$ and $-1$ (as 
$\hat{\alpha}_0>\max\{-P_0(1), P_0(-1)\}$) the result follows. 
\end{proof}

\noindent \textbf{Acknowledgment. } The second and third authors
are grateful to Institute of Mathematics at Academia Sinica in Taipei for
hospitality during their visit.


\end{document}